\documentclass[11pt,reqno]{amsart}

\usepackage[text={160mm,240mm},centering]{geometry}            
\usepackage{amssymb}
\usepackage{amsmath}
\usepackage{textcomp}
\usepackage{longtable}
\usepackage{float}
\usepackage{xcolor}
\newtheorem{theorem}{Theorem}[section]
\newtheorem{lemma}[theorem]{Lemma}
\newtheorem{proposition}[theorem]{Proposition}
\newtheorem{conjecture}[theorem]{Conjecture}
\newtheorem{corollary}[theorem]{Corollary}

\theoremstyle{definition}
\newtheorem{definition}[theorem]{Definition}

\theoremstyle{remark}
\newtheorem{remark}[theorem]{Remark}

\numberwithin{equation}{section}

\begin{document}

\title[unimodal  isolated complete intersection singularities]{Classification of unimodal  isolated complete intersection singularities  in positive characteristic}

\author{Hongrui Ma}
\address{Department of Mathematical Sciences,
Tsinghua University, Beijing, 100084, P. R. China.}
\email{mhr23@mails.tsinghua.edu.cn}

\author{Stephen S.-T. Yau}
\address{Hetao Institute of Mathematics and Interdisciplinary Sciences (HIMIS), Shenzhen 518000, Guangdong, P. R. China; Beijing Institute of Mathematical Sciences and Applications (BIMSA), Beijing 101408, P. R. China; Department of Mathematical Sciences, Tsinghua University, Beijing, 100084, P. R. China}
\email{yau@uic.edu}

\author{Huaiqing Zuo}
\address{Department of Mathematical Sciences,
Tsinghua University,
Beijing, 100084, P. R. China.}
\email{hqzuo@mail.tsinghua.edu.cn}
\thanks{Zuo is supported by NSFC Grant 12271280 and BJNSF Grant 1252009. Yau  is supported by the Tsinghua University Education Foundation.}





\begin{abstract}
In this paper we classify the unimodal isolated complete intersection singularities  in arbitrary characteristic under contact equivalence. The classification over $\mathbb{C}$ has already done by A. Dimca and C.G. Gibson. We continue and generalize their work. To complete the classification, we generalized the complete transversal method into positive characteristic field, which is also useful in many other classification problem.
\end{abstract}

\maketitle

\section{Introduction}
Classification is one of the oldest topic in singularity theory. The modality of singularities for real and complex hypersurfaces are first introduced by V.I. Arnold in \cite{arnold}. He also finished the classification of hypersurface singularities with small modality over $\mathbb{C}$ in \cite{arnold-class-C}. G.M. Greuel and H.D. Nguyen generalized the notion of modality to the algebraic setting in \cite{right-simple}, so that one can define modality over any algebraically closed field of arbitrary characteristic. They also classified simple (i.e. modality 0) hypersurface singularities in positive characteristic field under right equivalence. 

The classification of isolated complete intersection singularities (ICIS) under contact equivalence was studied in the 1980s. Assume $(X,0)=(f^{-1}(0),0)$ {\color{black}to be} a complete intersection germ with an isolated singularity defined by $f: F^n \rightarrow F^p$. We call such germs $I_{n,p}$. M. Giusti has shown that only $I_{1,1}$, $I_{2,2}$ and $I_{3,2}$ can be simple in \cite{germany-sing}. He then classified all simple ICIS over {\color{black}$F=\mathbb{C}$}. The classification of unimodal (i.e. modality 1) germs from plane to plane ($I_{2,2}$ case) was completed by A. Dimca and C.G. Gibson. C.T.C. Wall then classified the unimodal germs of $I_{n,p}$ with $n>p$. Recently, T.H. Pham, G. Pfister, and G.M. Greuel generalized the modality of hypersurface singularities to ICIS, and they classified zero-dimensional ICIS ($I_{2,2}$ case) over any algebraically closed field of arbitrary characteristic. In this paper, we continue their work. We use a different method to classify unimodal zero-dimensional ICIS of $I_{2,2}$ over any algebraically closed field of arbitrary characteristic.  In Corollary \ref{crit-of-unimodal}, we have shown that for all zero-dimensional ICIS, only $I_{2,2}$ and $I_{3,3}$ can be unimodal.  Unfortunately, the classification of $I_{3,3}$ seems very complicated and needs more new tools. 

In the previous work \cite{complete-unimodal-plane}, the finite determinacy theorem works well in zero characteristic fields. T.H. Pham and G.M. Greuel generalized it to positive characteristic in \cite{finitedeter}. Although the given boundary is sharp, it is  still insufficient to deal with complicated (i.e. higher order) problems. For example, let $h=(x^3+xy^3+y^5,x^2y+y^4+y^5)$ be the $5$-jet of an ICIS $f$ over an algebraically closed field $F$. If $char(F)=0$, then $h$ is $5$-determined. Therefore $f$ is contact equivalent to $(x^3+xy^3+y^5,x^2y+y^4+y^5)$, and then to $(x^3+xy^3+y^5,x^2y+y^4)$ after a transformation. Otherwise, we can only know that $h$ is $7$-determined in positive characteristic, which gives few idea about the form of $f$.

To solve the problem, we develop new tools. We generalized the complete transversal method introduced in \cite{complete-trans} into positive characteristic field, which is useful in many classification problem, see Corollary \ref{weight-ct-cor}. It can be used for semi-quasi homogeneous singularities over field of arbitrary characteristic. For the above example, $f$ is a semi-quasi homogeneous singularity with initial term $(x^3+xy^3,x^2y+y^4)$. {Using} our method we can show {that} $f$ is contact equivalent to $(x^3+xy^3+y^5,x^2y+y^4)$.

The main result is Theorem \ref{main-result}. Surprisingly, the classification result in positive characteristic turns out to be similar {as} the zero characteristic case {excepting} for some specical characteristic.

In the following, we set $F$ an algebraic closed field with arbitrary characteristic. $R=F[[x_1,\dots,x_n]]$, $\mathfrak{m}=\langle x_1,\dots,x_n \rangle \subset R$ and $p=char(F)$.

\section{Basic settings}

We first talk about some basic concepts of ICIS. Then we introduce the finite determinacy in positive characteristic in \cite{finitedeter}. After that, we generalize the complete transversal method into positive characteristic. The estimation of modality is also given in this section.

\subsection{Basic concepts}
\begin{definition}
    (1) An ideal $I \subset R$ defines a complete intersection if $I$ can be generated by $f_1,\dots,f_m$ with $f_i \in \mathfrak{m}$ for all $i$ such that $f_i$ is a non-zero divisor of $R/ \langle f_1,\dots,f_{i-1} \rangle$ for $i=1,\dots,m$. Then $\mathrm{dim}R/I=n-m.$

    (2) We call $f=(f_1,\dots,f_m)$ an isolated complete intersection singularity (ICIS) if $I=\langle f_1,\dots,f_m \rangle$ defines a complete intersection and  there exists $k \in \mathbb{N}$ such that $\mathfrak{m}^k \subset I+I_m(J(f))$, where $J(f)=(\frac{\partial f_i}{\partial x_j})_{ij}$ is the $m \times n$ Jacobian matrix and $I_m(J(f))$ is the {\color{black}ideal} generated by all $m \times m$ minors of $J(f)$. Denote $$I_{m,n}=\{f=(f_1,\dots,f_m) \in R^m \mid f \mathrm{\ is\ an \ ICIS\ with\ codimension\ }n-m \}.$$
    In this article, we {focus on} $I_{2,2},$ which denotes the zero-dimensional isolated complete intersection singularity in {the} plane.
\end{definition}

\begin{remark}
    Denote the Tjurina number $$\tau(f)=\mathrm{dim}_F R^m \bigg/ \bigg(\langle f_1,\dots,f_m \rangle \cdot R^m+ \langle \frac{\partial f}{\partial x_1},\dots,\frac{\partial f}{\partial x_n} \rangle \bigg).$$ Then a complete intersection $f$ is isolated if and only if $\tau(f) < \infty.$
\end{remark}

\begin{definition}
    The contact group $\mathcal{\mathcal{K}}$ is defined as $$\mathcal{\mathcal{K}}=GL(m,R) \rtimes Aut(R),$$ and the action of $\mathcal{\mathcal{K}}$ {on} $R^m$ is defined as $$(U,\phi,f) \mapsto U \cdot \phi (f),$$ with $U \in GL(m,R),\ \phi \in Aut(R),\ f=(f_1,\dots,f_m) \in R^m$ and $$\phi(f)=(f_1(\phi(\mathbf{x})),\dots,f_m(\phi(\mathbf{x})),$$ where {$\phi(\mathbf{x})=(\phi(x_1),\dots,\phi(x_n))$}.
\end{definition}

Let $f$ and $g$ define two isolated complete intersections {of the} same codimension $n-m$. $f$ is called contact equivalent to $g$, denoted by $f \sim g$, if $g \in \mathcal{K}f$, i.e. there exists $U \in GL(m,R)$ and $\phi \in Aut(R)$ such that $g=U\cdot \phi(f)$. 


\subsection{Tangent image and finite determinacy}
To classify the ICIS under contact equivalence, we need to work on jet spaces.
\begin{definition}\label{k-jet}
    (1) The $k$ -jet space of $R^m$ is defined as $J_k=R^m/\mathfrak{m}^{k+1}R^m$. For $f \in R^m$, the $k$-jet of $f$ is the image in $J_k$, denoted by $j_k(f)$. Let $\pi:J_l \rightarrow J_k$ be the natural projection and denote the kernel as $P_{k,l}.$ If $f \in J_k$ is a $k$-jet, we denote the submanifold $J_l(f)=f+P_{k,l}$.
   
    (2) We call $f$ is $k$-determined if for any $g \in R^m$ with $j_k(g)=j_k(f)$, we always have $g \sim f$.
\end{definition}

{Let $J_k=R^m/\mathfrak{m}^{k+1}R^m$ denote the $k$ -jet space of $R^m$.} Let $\mathcal{K}_k=\{(j_k(U),j_k(\phi)) \mid U \in GL(m,R),\ \phi \in Aut(R)\}$ be the $k$-jet algebraic group and the algebraic action of $\mathcal{K}_k$ on affine space $J_k$ is defined as $$(j_k(U),j_k(\phi),j_k(f)) \mapsto j_k(U\cdot \phi(f)).$$

The tangent space $T_e(\mathcal{K}_k)$ of algebraic group $\mathcal{K}_k$ has a natural Lie algebra structure (see \cite{actions-algebraic-group} {chapter} 4). The orbit map $\pi : \mathcal{K}_k \rightarrow \mathcal{K}_k \cdot f$ induces the tangent map $d\pi : Lie(\mathcal{K}_k) \rightarrow T_f(\mathcal{K}_k f)$. We denote the image of $d\pi$ as $\widetilde{T}_f(\mathcal{K}_k f)$, which {coinsides} with $Lie(\mathcal{K}_k) \cdot f$. When char$F=0$, $d\pi$ is surjective and 
    \begin{equation}\label{surj}
        \widetilde{T}_f(\mathcal{K}_k f)=T_f(\mathcal{K}_k f)
    \end{equation} (and therefore $\widetilde{T}_f(\mathcal{K} f)=T_f(\mathcal{K} f)$). But when char$F>0$, {(\ref{surj})} may not hold. For details one can see \cite{finitedeter}, section 2.

The tangent image is computed in \cite{finitedeter} as {follows}.

\begin{proposition}\label{tangent-image}
    The tangent image is identified with the submodule
    $$\widetilde{T}_f(\mathcal{K}_k f)=\bigg(\langle f_1,\dots,f_m \rangle \cdot R^m+\mathfrak{m} \cdot \langle \frac{\partial f}{\partial x_1},\dots,\frac{\partial f}{\partial x_n} \rangle+ \mathfrak{m}^{k+1}R^m \bigg) \bigg/ \mathfrak{m}^{k+1}R^m.$$ And the tangent image at $f$ to the orbit $\mathcal{K}f$ is the submodule
    $$\widetilde{T}_f(\mathcal{\mathcal{K}}f)=\langle f_1,\dots,f_m \rangle \cdot R^m+\mathfrak{m} \cdot \langle \frac{\partial f}{\partial x_1},\dots,\frac{\partial f}{\partial x_n} \rangle, $$ where $\langle f_1,\dots,f_m\rangle$ is regarded as an ideal of $R$, and $\langle \frac{\partial f}{\partial x_1},\dots,\frac{\partial f}{\partial x_n} \rangle$ is regarded as an ideal of $R^m$.
\end{proposition}

The finite determinacy {is} strongly related to tangent image.
\begin{theorem}\label{finite-determinacy}
    (cf. \cite{finitedeter} Theorem 3.2) Let $f=(f_1,\dots,f_m) \in R^m$. If there exists $k \in \mathbb{N}$ such that $$\mathfrak{m}^{k+2} \cdot R^m \subset \mathfrak{m} \cdot \widetilde{T}_f(\mathcal{K}f),$$ then $f$ is $(2k-\mathrm{ord}(f)+2)$-determined, where $\mathrm{ord}(f)=\mathrm{min}\{\mathrm{ord}(f_i) \mid i=1,\dots,m \}$. That is, for any $g \in R^m$ with $j_{2k-\mathrm{ord}(f)+2}(g)=j_{2k-\mathrm{ord}(f)+2}(f)$, we always have $g \sim f$.
\end{theorem}

\begin{remark}\label{finite-deter-tau}
    Let $f$ be an isolated complete intersection singularity, then $f$ is {($2\tau(f)-ord(f)+2$)-determined}.
\end{remark}

\subsection{Complete transversal}
In this part we introduce the complete transversal method.

Let $a=(a_1,\dots,a_n)$ be a given sequence of positive integers and $d=(d_1,\dots,d_m)$ a given sequence of non-negative integers. $f=(f_1,\dots,f_m)$ is said to be weighted homogeneous of degree $r$ (with respect to $(a;d)$) if $$f_i(t^{a_1}x_1,\dots,t^{a_n}x_n)=t^{r+d_i}f_i(x_1,\dots,x_n)$$ for any $t \in F$ and $i=1,2,\dots,m.$

Let $f$ be a $k$-jet in $J_k$ and weighted homogeneous of degree $0$ with respect to $(a;d)$. Moreover, assume  {
 \begin{equation}\label{condition-1.2}
        \mathrm{max}(d_i)<(k+1)\mathrm{min}(a_j) \mathrm{\ or\ }\mathrm{min}(d_i)>(k+1)\mathrm{max}(a_j) .
    \end{equation} }

For $l>k$, let $P_{k,l},\ J_l(f)$ be the {\color{black}subsets} of  $J_{l}$ defined in Definition \ref{k-jet}.(1). We have the following useful theorem from \cite{complete-unimodal-plane}.

\begin{theorem}\label{complete-cor}  
    For $f$ defined above, let $C \subset P_{k,l}$ be a linear subspace of $P_{k,l}$ satisfying $$P_{k,l} \subset C+\widetilde{T}_f(\mathcal{K}_{l}f)\cap P_{k,l},$$ we call $C$ a complete transversal. This complete transverse has the following property: every $g \in J_l(f)$ is in the same $\mathcal{K}_{l}$-orbit as some $l$-jet of the form $f+c$, for some $c \in C$.
\end{theorem}
\begin{proof}
    {We primarily refer to the proof of Proposition 1.3 in \cite{complete-unimodal-plane}.} Although the original proof is based on $\mathbb{C}$, nothing {\color{black}needs to be changed} in positive characteristic.
\end{proof}
\begin{remark}
    Condition {(\ref{condition-1.2})} is necessary. In \cite{InvariantsOH} Remark 1 there is a counterexample without {(\ref{condition-1.2})}: let $f=y^2+xy^3 \in F[[x,y]]$ with char$F=2$. One can check that $f$ {\color{black}does not satisfies} condition \ref{condition-1.2}. Computation shows $\widetilde{T}_f(\mathcal{K}_{l}f)=\mathfrak{m}^4/\mathfrak{m}^l$, which implies the complete transversal $C=0$ for any $l>4$ and then for any $g \in f+\mathfrak{m}^5$, $g \sim f$. But in fact we have $f+x^5 $ is not contact equivalent to $f$.
\end{remark}

\subsection{Complete transversal and homogeneous filtrations} In this part we introduce the generalization of the complete transversal method into weighted homogeneous filtrations in \cite{complete-trans}.

Let $F_{a,d}^rR^m$ denote the submodule of $R^m$ generated by the monomials of degree equal or greater than $r$ with respect to $(a;d)$. The sequence of $\{F_{a,d}^rR^m\}_{r \geq 0}$ defines a filtration of the module $F_{a,d}^0R^m$.

Next we {introduce} a filtration of contact group $\mathcal{K}$ compatible with the weighted filtration. For details one may see \cite{complete-trans} Section 2.3.

\begin{definition}\label{def-fil-group}
    (i) For $r \geq 0,$ define $$F^r\mathcal{R}=(I_n+F_{a,a}^rR^n) \cap \mathcal{R}.$$

    (ii) For $r \geq 0,$ define $$F^r\mathcal{C}=(I_{n+m}+F_{a\cup d,a\cup d}^r\widetilde{R}^{n+m}) \cap \mathcal{C},$$ where $$\widetilde{R}=F[[x_1,\dots,x_n,y_1,\dots,y_m]]$$ and $a\cup d$ denotes the $n+m$-tuple $(a_1,\dots,a_n,d_1,\dots,d_m).$

    (iii) Since the contact group $\mathcal{K}=\mathcal{R}\rtimes \mathcal{C},$ we define $$F^r\mathcal{K}=F^r\mathcal{R}\rtimes F^r\mathcal{C}.$$
\end{definition}
\begin{remark}
    For the survey of the standard Mather group $\mathcal{K},\mathcal{R},\mathcal{C}$, {one can refer to} \cite{mather-c-mapping}.
\end{remark}

\begin{proposition}\label{filtration-property}
    (i) {$F^r\mathcal{K}$ respects the filtration $\{F_{a,d}^r R^m\}$, i.e. for every $r,s\geq 0$, $(U,\phi) \in F^r\mathcal{K},\ f \in F^sR^m,\ U\cdot \phi(f) \in F^{s}R^m$, where $F^sR^m=F^s_{a,d}R^m$ is the submodule of $R^m$ generated by the monomials of degree equal or greater than $s$. }
    
    (ii) For $r,s,l \leq 0,$ the action of $F^r\mathcal{K}$ induces an action on $F^sR^m/F^{s+l}R^m$.

    (iii) The Lie {algebra} action satisfies the following: for any $f-g \in F_{a,d}^tR^m$ with $f,g \in F_{a,d}^0R^m$ and $l \in Lie(F_{a,d}^rR^m)$, we have $l \cdot f-l\cdot g \in F_{a,d}^{r+t}R^m.$
\end{proposition}

After a computation of tangent space (similar as the computation in \cite{finitedeter} Proposition 2.5), the tangent image of $F^r\mathcal{K}f$ can be regarded as

{
\begin{equation}\label{tangent-image-weight}
    \widetilde{T}_f(F_{a,d}^r\mathcal{K}\cdot f)= F_{a \cup d,d}^r(\langle f_1,\dots,f_m \rangle \cdot R^m)+\sum_j F_{a,a_j}^r(\mathfrak{m})\cdot \frac{\partial f}{\partial x_j} ,
\end{equation} }
{and clearly} $$\widetilde{T}_f(F_{a,d}^r\mathcal{K}\cdot f) \subset T_f(F_{a,d}^r\mathcal{K}\cdot f)$$ also holds. We denote $\langle f_1,\dots,f_m \rangle \cdot R^m+\langle \frac{\partial f}{\partial x_1},\dots,\frac{\partial f}{\partial x_n} \rangle$ as $\widetilde{T}_f^e(\mathcal{K}f)$ using the same sign in \cite{finitedeter}.

We have the similar complete transversal for $F^rR^m$ and $F^r\mathcal{K}$.

\begin{proposition}\label{weighted-fil-new1}
    Let $f \in F_{a,d}^0R^m$. $T$ is a subspace of $F_{a,d}^rR^m$ satisfying $$F_{a,d}^{k+1}R^m \subset T+Lie(F_{a,d}^1\mathcal{K})\cdot f+F_{a,d}^{k+2}R^m,$$then for any $g$ with $g-f \in F_{a,d}^{k+1}R^m$ is contact equivalent to $f+t+\overline{f}$ for some $t \in T$ and $\overline{f} \in F_{a,d}^{k+2}R^m$.
\end{proposition}

To prove the proposition, we first need a lemma, which can be seen as Taylor series in positive characteristic. We omit the proof.

\begin{lemma}\label{taylor}
    Let $f(x)=\sum_{\alpha}a_{\alpha}x^{\alpha} \in F[[x_1,\dots,x_n]]$ and $\xi=(\xi_1,\dots,\xi_n)\in F^n,$ then 
    \begin{equation}\label{taylor-eq}
            \begin{aligned}
                {f(x+\xi)}&=f(x)+\sum_{i=1}^n\xi_i \frac{\partial f}{\partial x_i}+\\ &terms\ of\ f(x)\ with\ more\ than\ two\ x_i\ replaced\ by\ \xi_i.
            \end{aligned}
    \end{equation}
\end{lemma}

Now we can prove Proposition \ref{weighted-fil-new1}.
\begin{proof}
Set $W=F_{a,d}^{k+1}R^m\backslash F_{a,d}^{k+2}R^m$ be the subspace of $F_{a,d}^{0}R^m$.

Claim: the following hold for any $f \in F_{a,d}^0R^m$ and $w \in W$:

    (i) $Lie(F^1\mathcal{K}) \cdot (f+w) -Lie(F^1\mathcal{K}) \cdot f \in F_{a,d}^{k+2}R^m$.

    (ii) $f+\{Lie(F^1\mathcal{K})\cdot f \cap W\}+F_{a,d}^{k+2}R^m \subset F^1\mathcal{K} \cdot f \cap \{f+W\}+F_{a,d}^{k+2}R^m.$ 

    {\color{black}\textit{Proof of claim}:}

    (i) follows from Proposition {\ref{filtration-property} (iii)}. For (ii), let $w=(w_1,\dots,w_m) \in \{Lie(F^1\mathcal{K})\cdot f \cap W\}$. Then $$w_j=\sum_{i=1}^n\xi_i \frac{\partial f_j}{\partial x_i}+\sum_{i=1}^m\lambda_{ji}f_i$$ is weighted homogeneous of degree $k+1$ with respect to $(a;d_j)$. Let $v_a(x)$ denote the minimal degree of monomials appeared in $x$ with respect to $(a;0)$, then $v_a(\xi_i)=a_i+k+1,\ v_a(\lambda_{ji}f_i)=d_j+k+1.$ Hence $v_a(\xi_i)\geq v_a(x_i)+k+1$ and $v_a(\lambda_{ji})\geq d_j-d_i+k+1.$

    Let $\phi \in Aut(R)$ such that $\phi(x_1,\dots,x_n)=(x_1-\xi_1,\dots,x_n-\xi_n).$ Set $U=(\lambda_{ji})\in M_{m\times m}(R).$ We have
    $$\phi(f+\sum_{i=1}^n\xi_i \frac{\partial f_j}{\partial x_i}-Uw)=f-\sum_{i=1}^n\xi_i \frac{\partial f_j}{\partial x_i}+h_1+\sum_{i=1}^n\xi_i \frac{\partial f_j}{\partial x_i}+h_2-Uw+h_3$$ by Lemma \ref{taylor}, where $h_i$ are the higher order terms defined in {(\ref{taylor-eq})}. We have $h_i \in F_{a,d}^{k+2}R^m$ for $i=1,2,3$ since $\xi_i$ appear twice in each terms and $v_a(\xi_i)\geq v_a(x_i)+k+1$. We also have $Uw\in F_{a,d}^{k+2}R^m$ since for each $j$, $v_a(\lambda_{ji}w_i)=v_a(\lambda_{ji})+v_a(w_i){\color{black}\geq} d_j-d_i+k+1+d_i+k+1\geq d_j+k+2$.

    Set $\widetilde{w}=h_1+h_2+h_3-Uw$. By above discussion, $\widetilde{w} \in F_{a,d}^{k+2}R^m.$ Without loss of generality, we can assume $v_a(f_1)\leq \dots \leq v_a(f_m).$ Then $U(0)$ is an upper triangular matrix whose principal diagonal element is zero. Hence $Id-U$ must be invertible. We have $f+w \sim (Id-U)\cdot \phi(f+w)=f+\widetilde{w}.$ Denote $g=(Id-U,\phi) \in \mathcal{K}$, we have $f+w=g^{-1}(f)+g^{-1}(\widetilde{w})$.

    To finish (ii), it remains to show $g \in F^1\mathcal{K}$. Since $\xi=(\xi_1,\dots,\xi_n) \in F_{a,a}^1R^n$ we have $\phi \in I_n+F_{a,a}^1R^n=F^1\mathcal{R}$ by Definition \ref{def-fil-group}. Similarly $Id-U \in I_m+F_{a\cup d,a \cup d}^1R^m \subset F^1\mathcal{C}.$ Now we get $g=(Id-U,\phi) \in F^1\mathcal{K}$. The claim is proved.

    We note that by claim (i) and (ii),
    \begin{equation}
        \begin{aligned}
            \bigcup_{t \in T}F^1\mathcal{K}\cdot (f+t+F_{a,d}^{k+2}R^m) &\supset \bigcup_{t \in T}\{f+t+Lie(F^1\mathcal{K})\cdot(f+t)\cap W+F_{a,d}^{k+2}R^m\}\\
            &=\bigcup_{t \in T}\{f+t+Lie(F^1\mathcal{K})\cdot f\cap W\}+F_{a,d}^{k+2}R^m\\
            &=f+T+Lie(F^1\mathcal{K})\cdot f\cap W+F_{a,d}^{k+2}R^m\\
            &=f+(T+Lie(F^1\mathcal{K})\cdot f) \cap W+F_{a,d}^{k+2}R^m\\
            &=f+W+F_{a,d}^{k+2}R^m.
        \end{aligned}
    \end{equation}

    That is, for any $g=f+w+\bar{g}$ with $w \in W$ and $\bar{g} \in F_{a,d}^{k+2}R^m$, $g$ is contact equivalent to $f+t+\bar{f}$ with $t \in T$ and $\bar{f} \in F_{a,d}^{k+2}R^m$.
\end{proof}

\begin{remark}
    In \cite{complete-trans}, Mather’s lemma implies (i) $\Rightarrow$ (ii). {However, the proof of Mather's lemma relies on analysis in the complex field.} We complete the proof in this special case without Mather's lemma. The proof refers M. Giusti's proof of \cite{germany-sing} Proposition 1.
\end{remark}

{Using} induction, we can show

\begin{proposition}\label{weight-ct}
    Suppose {that} $f$ is weighted homogeneous of weight $r$ with respect to $(a;d)$. Take $s>r$. Let $T$ be a subspace of $F_{a,d}^{r+1}R^m$ such that $$F_{a,d}^{r+1}R^m \subset T+\widetilde{T}_f(F_{a,d}^1\mathcal{K}\cdot f)+F_{a,d}^{s+1}R^m,$$ then any $g$ with $g-f \in F_{a,d}^{r+1}R^m$ is $F^1\mathcal{K}$-equivalent to $f+t+\phi$ where $t \in T$ and $\phi \in F_{a,d}^{s+1}R^m $.
\end{proposition}

\begin{proof}
    See Theorem 2.28 in \cite{complete-trans} and note that $$Lie(F_{a,d}^1\mathcal{K})\cdot f=\widetilde{T}_f(F_{a,d}^1\mathcal{K}\cdot f)$$ is shown above.
\end{proof}

For an isolated complete intersection singularity, $f$ is always finite determined. We can choose $s$ sufficiently large, then we have

\begin{corollary}\label{weight-ct-cor}
    Suppose {that} $f$ is an ICIS of weight $r$ with respect to $(a,d)$. Let $T$ be a subspace of $F_{a,d}^{r+1}R^m$ such that $$F_{a,d}^{r+1}R^m \subset T+\widetilde{T}_f(F_{a,d}^1\mathcal{K}\cdot f).$$ Then any $g$ with $g-f \in F_{a,d}^{r+1}R^m$ is contact equivalent to $f+t$ for some $t \in T$.
\end{corollary}

\subsection{Modality}

V.I. Arnol’d introduced the notion of modality in his famous \cite{arnold} as follows: The modality of a point $x \in X$ under the action of a Lie group $G$ on
a manifold $X$ is the smallest $m$ such that a sufficiently small neighborhood of $x$
may be covered by a finite number of orbit families of $m$ parameters. G.M. Greuel and H.D. Nguyen generalize the notion and give a detailed discuss in \cite{right-simple}, \cite{phdclassification}. For the definition of the modality of an ICIS, we refer to \cite{icissimple} Remark 1.13(3).

\begin{definition}
    An ICIS is called unimodal if $\mathcal{K}-mod(f)$, the $\mathcal{K}$-modality of $f$, is equal to 1.
\end{definition}

 In this section, we give some methods to estimate the lower bound of the modality and further give a criterion for non-unimodal.

First we give a lower bound by complete transversal in \cite{complete-unimodal-plane}, which is useful in next section.

Let $C$ be a complete transversal of $f$ in $J_l\ (l>k),$ for $a \in C$, we define
\begin{equation}
    \mathrm{cod}(f+a)=\mathrm{comdimension\ of\ }\widetilde{T}_f(\mathcal{K}_{l}f) \cap P_{k,l} \mathrm{\ in\ }P_{k,l}
\end{equation}
and
\begin{equation}
    \mathrm{cod}_0(f)=\mathrm{inf}_{a \in C}\{\mathrm{cod}(f+a)\}.
\end{equation}
Note that there exists a Zariski open subset $U \subset C$ such that $\mathrm{cod}(f+a)=\mathrm{cod}_0(f)$ if and only if $a \in U.$

\begin{proposition}\label{cod-modality}
    Let $f \in J_k$ be a $k$-jet of weighted homogeneous type with degree 0 with respect to $(a_1,\dots,a_n;d_1,\dots,d_m)$ and satisfies condition {(\ref{condition-1.2})}. Then for $a \in U,$ $f+a$ has modality $\mathrm{cod}_0(f)$ in $J_l(f)$ under the action of the subgroup $\mathcal{K}_l(f)$ of $\mathcal{K}_l$ which stabilize $f$. In particular, {any jet $h$ in $J_l(f)$ has $\mathcal{K}_l(f)-mod(h) \geq \mathrm{cod}_0(f)$ in $J_l$}.
\end{proposition}
\begin{proof}
    The idea mainly comes from \cite{complete-unimodal-plane} Proposition 1.4. We rewrite the proof using tangent image instead of tangent space for the sake of fields with positive characteristic.

    Find a subspace $\langle e_1,\dots,e_c \rangle \in P_{k,l}$ with $\langle e_1,\dots,e_c \rangle \oplus \widetilde{T}_f(\mathcal{K}_{l}(f+a)) \cap P_{k,l} = P_{k,l}$. Then $\mathrm{cod}(f+a)=c.$ Since $\widetilde{T}_f(\mathcal{K}_{l}(f+b))$ varies continuously for $b \in U$, we have $\langle e_1,\dots,e_c \rangle \cap \widetilde{T}_f(\mathcal{K}_{l}(f+b)) =\{0\}$ for $b$ in {a} Zariski open neighborhood $V$ of $a$.

    Consider 
    \begin{equation}
        \phi : (F^c,0) \rightarrow (J_l(f),f+a),\ (t_1,\dots,t_c) \mapsto f+a+\sum_{i=1}^c t_ie_i.
    \end{equation}
    We claim that: for any $\mathcal{K}_l(f)$-orbit $X$ in $J_l(f)$, $\phi^{-1}(X)$ consists of finitely many points in a neighbourhood of 0, hence $\phi$ is a minimal deformation of $f+a$ and $\mathcal{K}_l(f)-mod(f+a)=c=\mathrm{cod}_0(f).$

    For any $g \in J_l(f),$ write $g=f+\widetilde{g}=(f_1+\widetilde{g_1},\dots,f_m+\widetilde{g_m}),$ where $\widetilde{g}=(\widetilde{g_1},\dots,\widetilde{g_m})$ with weighted degree of $\widetilde{g_i}>d_i.$ We have 
    \begin{equation}
        \begin{aligned}
            g(x) \sim g_t(x)&=(t^{-d_1}g_1(t_1^{a_1}x_1,\dots,t_n^{a_n}x_n),\dots,t^{-d_n}g_m(t_1^{a_1}x_1,\dots,t_n^{a_n}x_n))\\
            &{=(f_1+t^{-d_1}\widetilde{g_1}(t_1^{a_1}x_1,\dots,t_n^{a_n}x_n),\dots,f_m+t^{-d_n}\widetilde{g_m}(t_1^{a_1}x_1,\dots,t_n^{a_n}x_n)).}
        \end{aligned}
    \end{equation}

    Then condition {(\ref{condition-1.2})} ensures any neighborhood of $f$ intersects $\mathcal{K}_l(f) \cdot g$, that is, any neighborhood of $f$ intersects all $\mathcal{K}_l(f)$-orbits in $J_l(f).$ Hence for any $\mathcal{K}_l(f)$-orbit $X$ in $J_l(f)$, there exists $b \in V,\ X=\mathcal{K}_l(f)\cdot (f+b).$ But we have {$\widetilde{T}_{f+b}(\mathcal{K}_l(f)\cdot (f+b))=\widetilde{T}_{f+b}(\mathcal{K}_l \cdot (f+b)) \cap P_{k,l},$} hence $\widetilde{T}_{f+b}(\mathcal{K}_l(f)\cdot (f+b)) \cap \langle e_1,\dots,e_c \rangle=\{0\}$, i.e. $\phi^{-1}(X)$ has only finitely many points in a neighborhood of 0. This finishes the proof.
\end{proof}

\begin{remark}\label{semicont-mod}
    In fact, we can show $\mathcal{K}-mod(f)$ is semicontinuous. That is, let $F(\mathbf{x},\mathbf{t}) \in F[\mathbf{t}][[\mathbf{x}]]$ such that $F_{\mathbf{t}_0} =F(\mathbf{x},\mathbf{t}_0)$ is an ICIS for a $\mathbf{t}_0 \in F^k$. Then there is a Zaraski open subset $U \in F^k$ such that $F(\mathbf{x},\mathbf{t})$ is an ICIS for any $\mathbf{t} \in U$. In addition, the sets $U_i=\{\mathbf{t} \in U \mid \mathcal{K}-mod(F(\mathbf{x},\mathbf{t})) \leq i\}$ are open for all $i \in \mathbb{N}$. In particular, let $mod_{min}=min\{\mathcal{K}-mod(F(\mathbf{x},\mathbf{t}))\mid \mathbf{t} \in F^k\}$, then $U_{min}=\{\mathbf{t} \in U \mid \mathcal{K}-mod(F(\mathbf{x},\mathbf{t}))=mod_{min}\}$ is open and dense.
\end{remark}

Next we use the following facts from \cite{phdclassification} to give a criterion for non-unimodal.

\begin{proposition}\label{modality}Let the algebraic group $G$ act on a variety $X$.

    {(1) If} the subvariety $X' \subset X$ is invariant under $G$ and $x \in X'$, then $$G-mod(x) \mathrm{\ in}\ X \geq G-mod(x) \mathrm{\ in}\ X'.$$

    {(2) Let} additionally the algebraic group $G'$ act on a variety $X'$ and let $p:X \rightarrow X'$ be a morphism of varieties. $p$ is open and $$ G\cdot x \subset p^{-1}(G' \cdot p(x)),\ \forall x \in X.$$ Then $$G-mod(x) \geq G'-mod(p(x)),\ \forall x \in X.$$
    
    {(3) If} $X$ is irreducible, for $x \in X$, we have $$G-mod(x) \geq \mathrm{dim}X-\mathrm{dim}G.$$
\end{proposition}

\begin{proposition}\label{icismod}
    Let $f \in I_{m,n}.$ Then $\mathcal{K}-mod(f)=\mathcal{K}_k-mod(j_k(f))$ for $k$ sufficiently large.
\end{proposition}

\begin{proof}
    See \cite{phdclassification}, Chapter 3.
\end{proof}

The following proposition is the main result of this section.
\begin{proposition}\label{crit-of-unimodal}
    Let $f \ in\ I_{n,n}$ with $\mathrm{ord}(f)=l$ and $f$ is unimodal. Then one of the following holds:\\
    (1) $n=2,\ l\leq 3$.\\
    (2) $n=3,\ l=2$.
\end{proposition}
\begin{proof} 
    Choose $k$ sufficiently large and let $X=\mathfrak{m}^l/\mathfrak{m}^{k+1}$. It follows from Proposition \ref{icismod} and \ref{modality}(1) that $$1=\mathcal{K}-mod(f)=\mathcal{K}_k-mod (f) \mathrm{\ in\ } J_k \geq \mathcal{K}_k-mod (f) \mathrm{\ in\ } X.$$

    Let $X'=\mathfrak{m}^l/\mathfrak{m}^{l+1}$. The action of $\mathcal{K}_k$ on $X$ induces the action of the algebraic group $\mathcal{K}'=GL(m,F)\times GL(m,F)$ on $X'$, and it can be easily checked that $p:X \rightarrow X'$ is open and $\mathcal{K}_k \cdot f \subset p^{-1}(\mathcal{K}' \cdot p(f))$. Then by Proposition \ref{modality}(2) we have $$\mathcal{K}_k-mod (f) \mathrm{\ in\ } X \geq \mathcal{K}'-mod (p(f)) \mathrm{\ in\ } X' .$$
    
    It is easy to see that $$\mathrm{dim}X'=n\binom{n-1+l}{l},$$ while for any $g \in X'$, $\mathrm{dim}(\mathcal{K}' \cdot g) \leq \mathrm{dim}\mathcal{K}'-1$, since $\{(a^l E_n,\frac{1}{a}E_n) \mid a \in F^{\times} \} \subset \mathcal{K}'$ stabilize $g$, {where $F^{\times}=F\backslash \{0\}$ denotes the units in $F$}. 
    
    After a small change of the proof of Proposition \ref{modality}(3) in \cite{phdclassification}, we have $$1 \geq \mathcal{K}'-mod (p(f)) \mathrm{\ in\ } X' \geq \mathrm{dim}X'-(\mathrm{dim}\mathcal{K}'-1),$$ which is $$1 \geq n\binom{n-1+l}{l}-(2n^2-1). $$
    The only solutions are $n=2,\ l \leq 3$ and $n=3,l=2$.
\end{proof}

\begin{corollary}\label{crit-of-unimodal}
    Let $f$ be {a} unimodal zero-dimensional ICIS, then $f \in I_{2,2}$ or $I_{3,3}$.
\end{corollary}

In the following we discuss the case $n=2$. The $n=3$ case is presented in the later article. From now on we assume $R=F[[x,y]]$, $\mathfrak{m}=\langle x,y \rangle$.

\section{The classification of order 2}
Some classification of order 2 ICIS has already been discussed in \cite{icissimple}, i.e. ICIS of modality 0. Here we continue their work and finish the classification of order 2 ICIS with modality 1. In this section $f=(f_1,f_2) \in R^2.$

First, we assume $charF \neq 2$. The following two {Propositions} are from \cite{icissimple}.

\begin{proposition}
    (i) If some of $j_2(f_i)$ is non-degenerate, then $f \sim (xy,x^n+y^m)$ for some $m,n \geq 2$, which is of modality 0.

    (ii) If $j_2(f)$ is degenerate, then
    \begin{equation}\label{ord2form}
        f \sim (x^2+\alpha y^s,\sum_{i \geq t}a_iy^i+x\sum_{j \geq q}b_jy^j),
    \end{equation}
    where $s \geq 3,\alpha \in \{0,1\},t\geq 2,q \geq 1$ and $ a_i,b_j \in F$.
\end{proposition}

\begin{proposition}\label{basic-trans}
    Let $f= ( f_{1}, f_{2}) = ( x^{2}+ \alpha y^{s}, \sum _{i\geq t}a_{i}y^{i}+ x\sum _{j\geq q}b_{j}y^{j})$ be an ICIS such that $s\geq 3$, $t\geq 2$, $q\geq 1$ and $\alpha \in \{ 0, 1\} .$
    
    (i) If $a_{i}= 0$ for all $i$ and $b_{q}\neq 0$, then $\alpha = 1$ and $f\sim ( x^{2}+ y^{s}, xy^{q})$.
    
    (ii) If $b_{j}=0$ for all $j$ and $a_t\neq0$, then $f\sim(x^2+\alpha y^s,y^t).$ If $\alpha=0$ or $t\leq s$, then $f \sim(x^2,y^t).$
    
    Assume now that $a_tb_q\neq 0.$
    
    (iii) If $t\leq q$, then $f\sim(x^2+\alpha y^s,y^t).$ If additionally $\alpha = 0$ or $t\leq s$ then $f \sim ( x^2, y^t) .$

    (iv) If $t>q$ and $\alpha=0$, then $f\sim(x^{2},y^{t}+xy^{q})$.
    
    (v) Let $t > q$ and $\alpha = 1$. Then $f \sim (x^{2}+ y^{s},y^{t}+exy^{q})$ for a suitable unit $e\in F[[y]]$. {If $2t - 2q - s \neq 0$ and $p = 0$ (or $p \nmid (2t - 2q - s)$), then $f\sim(x^{2}+y^{s},y^{t}+xy^{q})$}.
    
    (vi) If $ t=q+1$ and $p\nmid t$, then $f\sim(x^{2}+\alpha y^{s},y^{t})$.
\end{proposition}

Next we give a criterion of modality 1:

\begin{proposition}\label{ord2-modality1}
    Assume $f$ is of the form {(\ref{ord2form})}. Then if $s \geq 5, t \geq 6$ and $q \geq 3$, then $f$ is of modality at least 2.
\end{proposition}
\begin{proof}
    Let $g=j_3(f)=(x^2,0).$  Then an open subset in $J_5'(g)$ is formed by jets equivalent with $$h \sim (x^2+y^5, axy^3),$$ where{$J_5'(g)$} is formed by jets in $J_5(g)$ with $s \geq 5, t \geq 6$ and $q \geq 3$.

    
    For all $a,b \in F$, computation shows $(y^4,0),(0,y^4),(0,y^5) \notin \widetilde{T}_h(\mathcal{K}_5h)$. Hence the codimension of $\widetilde{T}_h(\mathcal{K}_5h) \cap P_{3,5}$ in $P_{3,5}$ is at least 2. By Proposition \ref{cod-modality}, $\mathcal{K}-mod(f) \geq \mathcal{K}_5-mod(f) \geq 2$.
\end{proof}
\begin{remark}
    By {(\ref{ord2form})}, $h \sim (x^2+\alpha y^5, y^6+axy^3+bxy^4+cxy^5)$. One can show, for example, applying $\phi(y)=y-\frac{c}{6}x$, that $h \sim (x^2+\alpha y^5, y^6+axy^3+bxy^4)$. This simplifies the computation of the codimension, but does not affect the result.
\end{remark}

The following proposition is also from \cite{icissimple}.
\begin{proposition}\label{simple-form}
    {Let $p=\text{char}(F)$.} The following ICIS are the only candidates for being simple ({i.e.} of modality 0):\\
    \textbf{(0) }$j_2(f_i)$ is non-degenerate, then $f \sim (xy,x^s+y^m),s,m \geq 2.$\\
    \textbf{(1) }$\forall a_i=0$, then $f \sim (x^2+y^3,xy^q)$, $q\geq 3$ or $f \sim (x^2+y^s,xy^2),\ s \geq 3$.\\
    \textbf{(2) }$\forall b_i=0,$\\
    \hspace*{1em}\textbf{(2.a)} $\alpha=0$ or $t \leq s$, then $f \sim (x^2,y^t),\ t=2,3,4$.\\
    \hspace*{1em}\textbf{(2.b) }$\alpha=1$ and $t>s$, then $f \sim (x^2+y^3,y^t),\ t \geq 4$.\\
    \textbf{(3)} $a_tb_q \neq 0$,\\
    \hspace*{1em}\textbf{(3.a)} $t\leq q,$\\
    \hspace*{1em}\hspace*{1em}\textbf{(3.a.i) }$\alpha=0$ or $t\leq s$, then $f \sim (x^2,y^t),t=2,3,4.$\\
    \hspace*{1em}\hspace*{1em}\textbf{(3.a.ii)} $\alpha=1$ and $t>s,$ then $f \sim (x^2+y^3,y^t),\ t \geq 4$.\\
    \hspace*{1em}\textbf{(3.b) }$t> q,$\\
    \hspace*{1em}\hspace*{1em}\textbf{(3.b.i)} $t=q+1,$\\
    \hspace*{1em}\hspace*{1em}\hspace*{1em}{\textbf{(3.b.i.1)} }$p \nmid t$,\\
    \hspace*{1em}\hspace*{1em}\hspace*{1em}\hspace*{1em}\textbf{(3.b.i.1.1) }$\alpha=0$ or $t \leq s$, then $f \sim (x^2,y^t),\ t=2,3,4$.\\
    \hspace*{1em}\hspace*{1em}\hspace*{1em}\hspace*{1em}\textbf{(3.b.i.1.2) }$\alpha=1$ and $t>s,$ then $f \sim (x^2+y^3,y^t),\ t \geq 4$.\\
    \hspace*{1em}\hspace*{1em}\hspace*{1em}{\textbf{(3.b.i.2)}} $p \mid t$,\\
    \hspace*{1em}\hspace*{1em}\hspace*{1em}\textbf{(3.b.i.2.1) }$\alpha=0$, then $f \sim (x^2,xy^2+y^3).$\\
    \hspace*{1em}\hspace*{1em}\hspace*{1em}\textbf{(3.b.i.2.2)} $\alpha=1$,\\
    \hspace*{1em}\hspace*{1em}\hspace*{1em}\hspace*{1em}\textbf{(3.b.i.2.2.1) }$s=3$, then $f \sim (x^2,xy^2+y^3)$ for $p=t=3$, or $f \sim (x^2+y^3,y^t+xy^{t-1})$ for $t \geq 4.$\\
    \hspace*{1em}\hspace*{1em}\hspace*{1em}\hspace*{1em}\textbf{(3.b.i.2.2.2)} $s>3,t=3,p=3$, then $f \sim (x^2,y^3+xy^2).$\\
    \hspace*{1em}\hspace*{1em}\textbf{(3.b.ii)} $t>q+1,$\\
    \hspace*{1em}\hspace*{1em}\hspace*{1em}\textbf{(3.b.ii.1)} $\alpha=0$,\\
    \hspace*{1em}\hspace*{1em}\hspace*{1em}\hspace*{1em}\textbf{(3.b.ii.1.1)} $q=1$, then $f \sim (xy,x^2+y^{2t-2}),\ t \geq 2.$\\
    \hspace*{1em}\hspace*{1em}\hspace*{1em}\hspace*{1em}\textbf{(3.b.ii.1.2)} $q=2$, then $f \sim (x^2+y^{2t-4},xy^2),\ t \geq 4$.\\
    \hspace*{1em}\hspace*{1em}\hspace*{1em}\textbf{(3.b.ii.2)} $\alpha=1$,\\
    \hspace*{1em}\hspace*{1em}\hspace*{1em}\hspace*{1em}\textbf{(3.b.ii.2.1)} $s \geq 3,\ q \leq 2,$ then $f \sim (xy,x^2+y^m)$ for some $m$ {and} $q=1$, and $f \sim (x^2+y^m,xy^2)$ for some $m$ {and} $q=2$.\\
    \hspace*{1em}\hspace*{1em}\hspace*{1em}\hspace*{1em}\textbf{(3.b.ii.2.2)} $s =3,\ t \geq q+3,$ then $f \sim (x^2+y^3,xy^q)$.\\
    \hspace*{1em}\hspace*{1em}\hspace*{1em}\hspace*{1em}\textbf{(3.b.ii.2.3)} {$s=3,\ t=q+2,$} then $f \sim (x^2+y^3,xy^q+y^{q+2})$.
\end{proposition}

Next we classify ICIS of modality 1 based on Proposition \ref{simple-form}.

\begin{proposition}\label{ord2-class-mod1}
    The following ICIS are the only candidates of modality 1:

\renewcommand\arraystretch{1.5}
\begin{longtable}{|c|c|c|}
\caption{}\label{table-ord2}
\\
\hline
Symbol&Form& condition  \\
\hline
$h_q$&$(x^2+y^4,xy^q)$&$q \geq 3$ \\
\hline
$i$&$(x^2,y^5)$& \\
\hline
$\widetilde{i}$&$(x^2,y^5+xy^3)$& \\
\hline
$i^5$&$(x^2,y^5+xy^4)$& $p=5$\\
\hline
$j_t$&$(x^2+y^4,y^t)$ &$t \geq 5$ \\
\hline
$\widetilde{j}_t$&$(x^2+y^4,y^t+xy^{t-1})$&$t \geq 5,\ p \mid t$ \\
\hline
$k_q$&$(x^2+y^4,y^{q+3}+xy^q)$&$q \geq 3$\\
\hline
$l_{q,\lambda}$&$(x^2+y^4,y^{q+2}+\lambda xy^q)$ & $q \geq 3, \lambda^2 \notin \{0, -1\}$\\
\hline
{$\widetilde{l}_{q,t,t'} $}&$(x^2+y^4,y^{q+2}+\lambda xy^q+uxy^t+xy^{t'}),\ where$ & $\lambda^2=-1, q \geq 3,$  \\
&$u=u_0+u_1y^p+u_2y^{2p}+\dots$& $ t \geq q+1, t' \geq t+1, $ \\
&& $p \mid t-q, p \nmid t'-q$\\
\hline
\end{longtable}
\end{proposition}

\begin{proof}
    \textbf{(0)} If $j_2(f_i)$ is non-degenerate, then $f \sim (xy,x^s+y^m),s,m \geq 2,$ which is simple.\\
    \textbf{(1)} $\forall a_i=0$, then $f \sim (x^2+y^4,xy^q)$, $q\geq 3$.\\
    \textbf{(2)} $\forall b_i=0,$\\
     \hspace*{1em}\textbf{(2.a)} $\alpha=0$ or $t \leq s$, then $f \sim (x^2,y^t),\ t=5$.\\
     \hspace*{1em}\textbf{(2.b)} $\alpha=1$ and $t>s$, then $f \sim (x^2+y^4,y^t),\ t \geq 5$.\\
    \textbf{ (3)} $a_tb_q \neq 0$,\\
    \hspace*{1em}\textbf{(3.a)} If $t\leq q,$ then $f \sim (x^2+\alpha y^s,y^t)$.\\
    \hspace*{1em}\hspace*{1em}\textbf{(3.a.i)} $\alpha=0$ or $t\leq s$, then $f \sim (x^2,y^t),t=5.$\\
    \hspace*{1em}\hspace*{1em}\textbf{(3.a.ii)} $\alpha=1$ and $t>s,$ then $f \sim (x^2+y^4,y^t),\ t \geq 5$.\\
    \hspace*{1em}\textbf{(3.b)} If $t> q,$\\
    \hspace*{1em}\hspace*{1em}\textbf{(3.b.i)} when $t=q+1,$\\
    \hspace*{1em}\hspace*{1em}\hspace*{1em}\textbf{(3.b.i.1)} If $p \nmid t$, then $f \sim (x^2+\alpha y^s,y^t)$ by Proposition \ref{basic-trans} (vi).\\
    \hspace*{1em}\hspace*{1em}\hspace*{1em}\hspace*{1em}\textbf{(3.b.i.1.1)} $\alpha=0$ or $t \leq s$, then $f \sim (x^2,y^t),\ t=5$.\\
    \hspace*{1em}\hspace*{1em}\hspace*{1em}\hspace*{1em}\textbf{(3.b.i.1.2)} $\alpha=1$ and $t>s,$ then $f \sim (x^2+y^3,y^t),\ t \geq 4$, which is simple.\\
    \hspace*{1em}\hspace*{1em}\hspace*{1em}\textbf{(3.b.i.2)} If $p \mid t$,\\
    \hspace*{1em}\hspace*{1em}\hspace*{1em}\hspace*{1em}\textbf{(3.b.i.2.1)} $\alpha=0$, then $f \sim (x^2,y^3+xy^2)$, which is of modality 0 by Proposition \ref{simple-form}.\\
    \hspace*{1em}\hspace*{1em}\hspace*{1em}\hspace*{1em}\textbf{(3.b.i.2.2)} $\alpha=1$, then $f \sim (x^{2}+ y^{s},y^{t}+exy^{q})$ by Proposition \ref{basic-trans} (v). And we have $s=3,4$ or $p=t=3,5$ by Proposition \ref{ord2-modality1}.\\
    \hspace*{1em}\hspace*{1em}\hspace*{1em}\hspace*{1em}\hspace*{1em}\textbf{(3.b.i.2.2.1)} If $s=3$, then $f \sim (x^2,xy^2+y^3)$ for $p=t=3$, or $f \sim (x^2+y^3,y^t+xy^{t-1})$ for $t \geq 4$ by Proposition \ref{basic-trans} (v), which is of modality 0.\\
    \hspace*{1em}\hspace*{1em}\hspace*{1em}\hspace*{1em}\hspace*{1em}\textbf{(3.b.i.2.2.2)} If $s=4$, then $f \sim (x^2+y^4,y^t+xy^{t-1})$ for $t \geq 4$ by Proposition \ref{basic-trans} (v). If $p=t=3,$ then $f \sim (x^2+y^4-y(y^3+xy^2),y^3+xy^2) \sim (x^2,y^3+xy^2)$ additionally, which is simple as shown in Proposition \ref{simple-form}. Hence $f \sim (x^2+y^4,y^t+xy^{t-1})$ for $p \mid t, t\geq 5$.\\
    \hspace*{1em}\hspace*{1em}\hspace*{1em}\hspace*{1em}\hspace*{1em}\textbf{(3.b.i.2.2.3)} If $p=t=3,\ s>4,$ same as the process in \cite{icissimple} Proposition 2.5, $f \sim (x^2,y^3+xy^2)$, which is simple.\\
    \hspace*{1em}\hspace*{1em}\hspace*{1em}\hspace*{1em}\hspace*{1em}\textbf{(3.b.i.2.2.4)} If $p=t=5,\ s>4$, then 
    {\begin{equation}
        \begin{aligned}
            f &\sim (x^{2}+ y^{s},y^{5}+e(y)xy^{4}),\ where\ e(y)\in F[[x,y]]\ is\ a\ unit\\
             &\sim (x^{2}+ e_0y^{s},y^{5}+xy^{4}),\ e_0\in F \\
             &\sim (x^2+e_0y^s-e_0y^{s-5}(y^5+xy^4),y^5+xy^4) \\
             &\sim (x^2-e_0xy^{s-1},y^5+xy^4) \\
             &\sim ((x-\frac{e_0}{2}y^{s-1})^2-\frac{e_0^2}{4}y^{2s-2},y^5+xy^4).
        \end{aligned}
    \end{equation}}
    Using the automorphism $\phi(x)=x-\frac{e_0}{2}y^{s-1},\phi(y)=y$, then 
    \begin{equation}
        \begin{aligned}
            f \sim (x^2-\frac{e_0^2}{4}y^{2s-1},y^5+xy^4+\frac{e_0}{2}y^{s+3}) \\
            \sim (x^2-\frac{e_0^2}{4}y^{2s-1},y^5+\widetilde{e}xy^4),
        \end{aligned}
    \end{equation}
    where $\widetilde{e}=\frac{1}{1+\frac{e_0}{2}y^{s-2}}$. Applying $\phi(x)=\frac{1}{\widetilde{e}}x$ and $\phi(y)=y$ , we have $f \sim (x^2+e_1y^{2s-1},y^5+xy^4),$ {where $e_1=-\frac{\widetilde{e}^2e_0^2}{4}$.} Repeating the process, we get $f \sim (x^2,y^5+xy^4)$ with $p=char F=5$.\\
    \hspace*{1em}\hspace*{1em}\textbf{(3.b.ii)} Now assume $t>q+1,$\\
    \hspace*{1em}\hspace*{1em}\hspace*{1em}\textbf{(3.b.ii.1)} If $\alpha=0$, by Proposition \ref{ord2-modality1}, we have {$1 \leq q \leq 2$} or $q=3,t=5$.\\
    \hspace*{1em}\hspace*{1em}\hspace*{1em}\hspace*{1em}\textbf{(3.b.ii.1.1)} If $q=1$, then $j_2(f_2)$ is non-degenerate, hence $f$ is simple.\\
    \hspace*{1em}\hspace*{1em}\hspace*{1em}\hspace*{1em}\textbf{(3.b.ii.1.2)} If $q=2$, then $f \sim (x^2+y^{2t-4},xy^2),\ t \geq 4$ same as Propositon \ref{simple-form}, which is also simple.\\
    \hspace*{1em}\hspace*{1em}\hspace*{1em}\hspace*{1em}\textbf{(3.b.ii.1.3)} If $q=3,t=5$, then $f \sim (x^2,y^5+xy^3)$.\\
    \hspace*{1em}\hspace*{1em}\hspace*{1em}\textbf{(3.b.ii.2)} If $\alpha=1$, we have $f \sim (x^2+y^s,y^t+exy^q)$ by Proposition \ref{basic-trans} (v). By Proposition \ref{ord2-modality1}, we have $s=3,4$ or $t=2,3,4,5$ or $q=1,2$.\\
    \hspace*{1em}\hspace*{1em}\hspace*{1em}\hspace*{1em}\textbf{(3.b.ii.2.1)} If $q=1,2$ holds, then $f$ is simple as shown in Proposition \ref{simple-form}.\\
    \hspace*{1em}\hspace*{1em}\hspace*{1em}\hspace*{1em}\textbf{(3.b.ii.2.2)} If $s=3$, then $f \sim (x^2+y^3,xy^q)$ or $(x^2+y^3,y^{q+2}+xy^q)$ {for $q \geq 3$}, which is also simple. \\
    {\color{black}\hspace*{1em}\hspace*{1em}\hspace*{1em}\hspace*{1em}\textbf{(3.b.ii.2.3)}} If $s=4,q>2$,\\
    \hspace*{1em}\hspace*{1em}\hspace*{1em}\hspace*{1em}\hspace*{1em}\textbf{(3.b.ii.2.3.1) }If $t \geq q+4$, then $f \sim (x^2+y^4,y^t+e(y)xy^q-y^{t-4}(x^2+y^4))=(x^2+y^4,(e(y)-xy^{t-q-4})xy^q)\sim (x^2+y^4,xy^q),$ {where $e(y)=b_q+b_{q+1}y+\dots$ is a unit in $F[[x,y]]$ and $b_q,b_{q+1}$ are defined in \eqref{ord2form}.}\\
    \hspace*{1em}\hspace*{1em}\hspace*{1em}\hspace*{1em}\hspace*{1em}\textbf{(3.b.ii.2.3.2)} If $t = q+3$, then $f \sim (x^2+y^4,y^{q+3}+xy^q)$ by Proposition \ref{basic-trans} (v) since in this case $2t-2q-s=2$.\\
    \hspace*{1em}\hspace*{1em}\hspace*{1em}\hspace*{1em}\hspace*{1em}\textbf{(3.b.ii.2.3.3)} If $t=q+2$. Let $g=(x^2+y^4,y^{q+2}+\lambda xy^q)$, where $\lambda=e(0) \in F^{\times}$, {$e(0)$ is the constant term of $e(y)$ in \textbf{(3.b.ii.2.3.1)}.} Then $g$ is weighted homogeneous of degree 0 with respect to $(a;d)$, where $a=(2,1),d=(4,q+2)$.\\
    \hspace*{1em}\hspace*{1em}\hspace*{1em}\hspace*{1em}\hspace*{1em}\hspace*{1em}\textbf{(3.b.ii.2.3.3.1)} If $\lambda^2 \neq -1$, we claim that $F_{a,d}^1R^2 \subset \widetilde{T}_g(F^1\mathcal{K}g)$ (the proof will be presented later). Hence by Proposition \ref{weight-ct}, $f \sim (x^2+y^4, y^{q+2}+\lambda xy^q),\ \lambda \in F^{\times}$ and $ \lambda^2+1 \neq 0.$\\\textbf{}
    \hspace*{1em}\hspace*{1em}\hspace*{1em}\hspace*{1em}\hspace*{1em}\hspace*{1em}\textbf{(3.b.ii.2.3.3.2)} If $\lambda^2+1=0$, from the proof of (3.b.ii.2.3.3.1), we can see $$F_{a,d}^1R^2 \subset T+\widetilde{T}_g(F^1\mathcal{K}g)$$ with $T=span \langle (0,xy^t) \mid t>q \rangle$. Then by Proposition \ref{weight-ct}, $$f \sim (x^2+y^4,y^{q+2}+\lambda xy^q+e(y)xy^t)$$ for {$\lambda^2+1=0,t>q$} and $e(y)$ is a unit in $F[[y]]$. In fact, we can show that for some $e(y)$ and $l$, $f \sim (x^2+y^4,y^{q+2}+\lambda xy^q+xy^l)$, while for others, we get a family of unimodal ICIS with the form
    \begin{equation}\label{family-x2y4}
        (x^2+y^4,y^{q+2}+\lambda xy^q+u(y)xy^t+xy^{t'}),
    \end{equation}
    where $u(y)=u_0+u_1y^p+u_2y^{2p}+\dots$ is a unit, $t'>t \geq q+1,\ p \mid t-q,\ p \nmid t'-q$. The details will be shown later.
    \\
    \hspace*{1em}\hspace*{1em}\hspace*{1em}\hspace*{1em}\textbf{(3.b.ii.2.4)} The last case remained is $t=5,q=3,s>4$. Then $f \sim (x^2+y^s,y^5+exy^3)$. Use the method same as \textbf{(3.b.i.2.2.4)}, one can show $f \sim (x^2,y^5+xy^3)$. 
\end{proof}

$Proof\ of\ the\ claim\ in\ $\textbf{(3.b.ii.2.3.3.1)}:

We need to show $F_{a,d}^1R^2 \subset \widetilde{T}_g(F^1\mathcal{K}g)$ for $g=(x^2+y^4,y^{q+2}+\lambda xy^q)$, where $a=(2,1),d=(4,q+2)$ and $\lambda^2 \neq -1.$ 

Denote {$e_1=(x^2+y^4,0),e_2=(0,x^2+y^4),e_3=(y^{q+2}+\lambda xy^q,0),e_4=(0,y^{q+2}+\lambda xy^q),e_5=(2x,\lambda y^q),e_6=(4y^3,(q+2)y^{q+1}+\lambda q xy^{q-1})$. }

Set the weight $weight(x)=2,weight(y)=1,weight(x^2+y^4)=4,weight(y^{q+2}+\lambda xy^q)=q+2.$ Then by {(\ref{tangent-image-weight})}, $\widetilde{T}_g(F^1\mathcal{K}g)$ has the element $x^iy^je_k$ with:\\
(a)  {$k=1,4,\ weight(x^iy^je_k) \geq (5,q+3)$}\\
(b) $k=5,\ weight(x^iy^j) \geq 3$ \\
(c) $k=6,\ weight(x^iy^j) \geq 2$.

We have $$xe_6=(4xy^3,(q+2)xy^{q+1}+\lambda q x^2y^{q-1}) \in \widetilde{T}_g(F^1\mathcal{K}g)$$ and $$ye_4-y^{q-1}e_2 =(0,-x^2y^{q-1}+\lambda xy^{q+1}) \in \widetilde{T}_g(F^1\mathcal{K}g).$$ Therefore $$xe_6+\lambda q (ye_4-y^{q-1}e_2)=(4xy^3,(\lambda^2q+q+2)xy^{q+1})\in \widetilde{T}_g(F^1\mathcal{K}g).$$

We also have $$y^{q-1}e_2=(0,x^2y^{q-1}+y^{q+3}) \in \widetilde{T}_g(F^1\mathcal{K}g),$$ hence, $$xe_6-\lambda q y^{q-1}e_2=(4xy^3,(q+2)xu^{q+1}-\lambda q y^{q+3}) \in \widetilde{T}_g(F^1\mathcal{K}g).$$ Then $$(q+2)(xe_6-\lambda q y^{q-1}e_2)+\lambda q y^2e_6=(4(q+2)xy^3+4\lambda q y^5,((q+2)^2+\lambda^2q^2)xy^{q+1}) \in \widetilde{T}_g(F^1\mathcal{K}g).$$ Denoted it by $e_7$.

Note that $$xye_5-2ye_1=(-2y^5,\lambda x y^{q+1}),\ y^3e_5-\lambda y e_4=(2xy^3,-\lambda^2 xy^{q+1}) \in \widetilde{T}_g(F^1\mathcal{K}g),$$ we have
\begin{equation}
    \begin{aligned}
        &e_7+ 2\lambda q (xye_5-2ye_1)-2(q+2)(y^3e_5-\lambda y e_4) \\
          =&(0,((q+2)^2+\lambda^2q^2+2\lambda^2q+2\lambda^2(q+2))xy^{q+1}) \\
          =&(0,(q+2)^2(\lambda^2+1)xy^{q+1}) \in \widetilde{T}_g(F^1\mathcal{K}g).
    \end{aligned}
\end{equation}

Since $\lambda^2+1 \neq 0$, we have $(0,xy^{q+1}) \in \widetilde{T}_g(F^1\mathcal{K}g)$. The other elements in $F_{a,d}^1R^2$ follow easily.

$Proof\ of\ the\ claim\ in\ $\textbf{(3.b.ii.2.3.3.2)}:

We have $f \sim (x^2+y^4, y^{q+2}+\lambda xy^q+e(y)xy^t)$. If $p \nmid t-q$, we can use an `$\alpha,\beta$-trick' based on the Implicit Function Theorem to show that $f \sim (x^2+y^4, y^{q+2}+\lambda xy^q+xy^t)$ for $t \geq q+1$. See Remark \ref{alpha-beta}. For the case $p \mid t-q$, {we write $$e(y)=\sum_{i \geq 0}e_iy^{q+ip}+\sum_{j\geq 0,\ p\nmid t'+j}e'_jy^{t'+j},$$ then $f \sim (x^2+y^4, y^{q+2}+\lambda xy^q(1+e_1y^p+e_2y^{2p}+\dots)+xy^{t'}(e_0'+e_1'y+e_2'y^2+\dots))$ with $p \nmid t'-q$.} Therefore we can use $\alpha,\beta$-trick again to reduce $e_0'+e_1'y+\dots$. Then we have $f \sim (x^2+y^4,y^{q+2}+\lambda xy^q+u(y)xy^t+xy^{t'})$, where $u(y)=u_0+u_1y^p+u_2y^{2p}+\dots$ as we want.

\section{The classification of order 2 in $charF= 2$}

In this section, we will show:

\begin{proposition}\label{ord2-class-char2}
    A unimodal ICIS of order 2 in any field with characteristic equal to $2$ must have the form in Table \ref{table-ord2-char2}.

\renewcommand\arraystretch{1.5}
\begin{longtable}{|c|c|c|}
\caption{}\label{table-ord2-char2}
\\
\hline
Symbol&Form& condition  \\
\hline
$h^2_{\lambda}$&$(x^2+\lambda xy^2,y^3)$& $\lambda \in \{0,1\}$ \\
\hline
$i^2_k$&$(x^2+y^k,xy^2)$&$k \geq 3, k$ is odd \\
\hline
$i^2_{k,\lambda}$&$(x^2+y^k+\lambda y^{k+1},xy^2)$& $k \geq 3, k$ is even, $\lambda \in \{0,1\}$ \\

\hline
$j^2_{\lambda}$ & $(x^2+y^3,y^4+\lambda xy^3)$ & $\lambda \in \{0,1\}$\\
\hline
$k^2_{\lambda,\mu}$ &$(x^2+\lambda xy^3,y^4+\mu xy^3)$ & $\lambda, \mu \in \{0,1\}$\\
\hline
$l^2$ & $(x^2+xy^2,y^4)$ &  \\
\hline
$m^2_{s}$ & $(x^2+y^s,xy^3)$ & $s \geq 3,\ s$ is odd \\
\hline
$m^2_{s,\lambda}$ & $(x^2+\lambda x^2y+ y^s,xy^3)$ & $\lambda \in F,\ s \geq 4,\ s$ is even \\
\hline
$n^2_{s}$ & $(x^2+xy^2+y^s,xy^3)$ & $s \geq 3,\ s$ is odd\\
\hline
$\widetilde{n}^2_{s,\lambda}$ & $(x^2+xy^2+\lambda x^2y+y^s,xy^3)$ & $\lambda \in F,\ s \geq 4,\ s$ is even \\
\hline

\end{longtable}

\end{proposition}

The following result is from \cite{icissimple}:

\begin{proposition}\label{char2-ord2-simple}
    Let $charF=2,$ $f \in I_{2,2}$ and $ord(f)=2$. Then one of the following cases occurs:\\
    (a) $f \sim (xy,g)$ for some {$g \in \mathfrak{m}^2$}. In this case, $f \sim (xy,x^m+y^n)$ for some $m,n \geq 2$, which is simple.\\
    (b) $f \sim (x^2+h,g)$ for $h \in \mathfrak{m}^3$ and $g \in \mathfrak{m}^2$. Moreover, if $g \notin \mathfrak{m}^3$, then $f$ is simple.
\end{proposition}

Moreover, similar {to} Proposition \ref{ord2-modality1}, we can show
\begin{proposition}\label{ord2-mod1-char2}
    If $g \in \mathfrak{m}^5$, then the modality of $f \sim (x^2+h,g)$ is at least $2$.
\end{proposition}
\begin{proof}
    For $g \in \mathfrak{m}^5$, then $l=j_4(f)$ is of the form
    $$l \sim (x^2+axy^2+by^3+cxy^3+dy^4,0) \sim (x^2+axy^2+by^3+dy^4,0).$$
    It is easy to compute that the codimension of $\widetilde{T}_l(\mathcal{K}_4l)$ in $P_{2,4}$ is at least $2$ ($(0,xy^3),(0,y^4) \notin P_{2,4}$). Therefore $\mathcal{K}-mod(f) \geq \mathcal{K}_4-mod(l) \geq 2$ by Proposition \ref{cod-modality}.
\end{proof}

Therefore, we need to work on the case $f \sim (x^2+h,g)$ with $g \in \mathfrak{m}^4 \backslash \mathfrak{m}^5$.

\begin{proposition}
    If $g \in \mathfrak{m}^3 \backslash \mathfrak{m}^4$, then $f \sim (x^2+\lambda xy^2,y^3),\ \lambda \in \{0,1\}$ or $f \sim (x^2+y^k,xy^2),\ k$ is odd or $f \sim (x^2+y^k+\mu y^{k+1},xy^2),\ k$ is even, $\mu \in \{0,1\}$.
\end{proposition}
\begin{proof}
    Set $j_3(g)=ax^3+bx^2y+cxy^2+dy^3$.

    If $d \neq 0$, let $f_0=(x^2,dy^3) \sim (x^2,y^3)$ be the weighted $0$-jet of $f$ {\color{black}with respect to} $(a,d)$, where $a=(4,3),\ d=(8,9)$. Then $$F_{a,d}^1R^2 \subset span\langle(xy^2,0) \rangle+\widetilde{T}_{f_0}(F^1\mathcal{K}f_0).$$ By Proposition \ref{weight-ct}, $f \sim (x^2+cxy^2,y^3) \sim (x^2+\lambda xy^2,y^3),\ \lambda \in \{0,1\}$.

    If $d=0$ and $c \neq 0$, we still choose $f_0=(x^2,cxy^2) \sim (x^2,xy^2)$ be the weighted $0$-jet of $f$ with respect to $(a,d)$, where $a=(4,3),\ d=(8,10)$. In this case $$F_{a,d}^1R^2 \subset span\langle(y^k,0),\ k\geq 4 \rangle+\widetilde{T}_{f_0}(F^1\mathcal{K}f_0).$$ Then $f \sim (x^2+e(y)y^k,xy^2) \sim (e(y)^{-1}x^2+y^k,xy^2),\ e(y) \in F[[y]]$ is a unit. 

    If $k$ is odd, then there exists $\widetilde{e}(y)^k=e(y)$. Apply $\phi(x)=\widetilde{e}(y)x,\ \phi(y)=y$, then $f \sim (x^2+y^k,xy^2)$.

    If $k$ is even, write $e(y)^{-1}=e_0+e_1y+\dots,$ then 
    \begin{equation}
        \begin{aligned}
            f &\sim (e_0x^2+e_1x^2y+e_2x^2y^2+\dots+y^k,xy^2) \\
            &\sim (e_0x^2+e_1x^2y+y^k,xy^2)\\
            &\sim ((e_0x^2+e_1x^2y+y^k)(1-\frac{e_1}{e_0}y),xy^2)\\
            &\sim (e_0x^2+y^k-\frac{e_1}{e_0}y^{k+1},xy^2)\\
            &\sim (x^2+y^k+\lambda y^{k+1},xy^2),
        \end{aligned}
    \end{equation}
    where $\lambda \in \{0,1\}$.

    If $c=d=0$, then $j_3(f) \sim (x^2+h,0)$. Then $g \in \mathfrak{m}^4$, a contradiction. 
\end{proof}

\begin{proposition}
    If {$g \in \mathfrak{m}^4$}, then $f$ is equivalent to other forms in Table \ref{table-ord2-char2}, {\color{black}that is, $j_\lambda^2$, $k_{\lambda,\mu}^2$, $l^2$, $m_s^2$, $m_{s,e_1}^2$, $n_s^2$ or $\widetilde{n}^2_{s,e_1}$.}
\end{proposition}

\begin{proof}
    Computing the complete transversal of $j_2(f)=(x^2,0)$, we have $f \sim (x^2+a(y)xy^r+b(y)y^s,c(y)xy^u+d(y)y^v)$ with $a(y),b(y),c(y),d(y)$ are units or $0$. If $a(y)$ (resp. $b(y),c(y),d(y)$)$=0$, we regard as $r$ (resp. $s,u,v$) $= \infty$. By Proposition \ref{ord2-mod1-char2}, we have either $u=3$ or $v=4$. 

    (i) $v=4,s=3$. Then $f \sim (x^2+a(y)xy^r+b(y)y^3,c(y)xy^u+d(y)y^4)$. Take $l=(x^2+b_0y^3,d_0y^4) \sim (x^2+y^3,y^4)$ be the weighted 0-jet with respect to $a=(3,2),d=(6,8)$. We have $F_{a,d}^1 \subset span\langle (0,xy^3) \rangle +\widetilde{T}_{l}(F^1\mathcal{K}l)$. {Therefore $f \sim (x^2+y^3,y^4+\lambda xy^3)$ for $\lambda \in \{0,1\}$, which is $j_\lambda^2$ in Table \ref{table-ord2-char2}.} 

    (ii) $v=4,s>3$. If $r \geq 3$, then we choose $l=(x^2,y^4)$ be the weighted 0-jet of $f$ with respect to $a=(2,1),d=(4,4)$. We have $F_{a,d}^1 \subset span\langle (xy^3,0),(0,xy^3) \rangle +\widetilde{T}_{l}(F^1\mathcal{K}l)$. Thus $f \sim (x^2+\lambda xy^3,y^4+\mu xy^3)$, $\lambda,\mu \in \{0,1\}$ after a scaling. {That is, $f \sim k^2_{\lambda,\mu}$ in Table \ref{table-ord2-char2}}. If $r=2$, then the weighted 0-jet of $f$ with respect to $a=(2,1),d=(4,4)$ becomes $l=(x^2+a_0xy^2,y^4) \sim (x^2+xy^2,y^4)$. Computation shows that $F_{a,d}^1 \subset \widetilde{T}_{l}(F^1\mathcal{K}l)$, {then we have $f \sim (x^2+xy^2,y^4)\sim l^2$, which is in Table \ref{table-ord2-char2}.}

    (iii) $v>4$. Then we have $s=3$ and then $f \sim (x^2+a(y)xy^r+b(y)y^s,xy^3+d(y)y^v)$. Choose $l=(x^2,xy^2)$ be the weighted 0-jet of $f$ with respect to $a=(1,1),d=(2,4)$. We have $F_{a,d}^1 \subset span\langle (xy^2,0),(y^k,0)|\ k \geq 3 \rangle +\widetilde{T}_{l}(F^1\mathcal{K}l)$. Then $f \sim (x^2+\mu xy^2+e(y)y^s,xy^3)$, $\mu \in \{0,1\}, s \geq 3$. After a scaling we can furthermore assume $e(0)=1$.

    If $\mu=0$ and $s$ is odd, using $\alpha,\beta$-trick in Remark \ref{alpha-beta}, we have $f \sim (x^2+y^s,xy^3)\sim m^2_s$. If $\mu=0$ and $s $ is even, we have  {
    \begin{equation}
        \begin{aligned}
            f &\sim (e(y)^{-1}x^2+y^s,xy^3) \\
            &\sim ((1+e_1y+e_2y^2+\dots)x^2+y^s,xy^3)\\
            &\sim (x^2+e_1x^2y+e_2x^2y^2+y^s,xy^3),
        \end{aligned} 
    \end{equation} where $e(y)^{-1}=e_0+e_1y+e_2y^2+\dots,$ and $e_3$-term vanishes since $e_3x^2y^3$ is killed by the $xy^3$ term in the second component.}
    We apply $$\phi(x)=\frac{x}{1+e_2^{\frac{1}{2}}y},\phi(y)=y,$$ then $f \sim (x^2+e_1(1+e_2y^2)^{-1}x^2y+y^s,xy^3)\sim (x^2+e_1x^2y+y^s,xy^3)$, $e_1 \in F$. That is, $f \sim m^2_{s,e_1}$

    If $\mu=1$ and $s$ is odd, we have $f \sim (x^2+xy^2+y^s,xy^3)\sim n^2_s$ by $\alpha,\beta$-trick. If $\mu=1$ and $s$ is even, as above, we have
    \begin{equation}
        \begin{aligned}
            f &\sim (e(y)^{-1}(x^2+xy^2)+y^s,xy^3) \\
            &\sim ((1+e_1y+e_2y^2+\dots)x^2+y^s,xy^3)\\
            &\sim (x^2+xy^2+e_1x^2y+e_2x^2y^2+y^s,xy^3).
        \end{aligned}
    \end{equation}
    Then apply $$\phi(x)=\frac{x}{1+e_2^{\frac{1}{2}}y},\phi(y)=y,$$ we have 
    \begin{equation}
        \begin{aligned}
            f &\sim (x^2+(1+e_2^{\frac{1}{2}}y)^{-1}xy^2+e_1(1+e_2y^2)^{-1}x^2y+y^s,xy^3)\\
            &\sim (x^2+xy^2+e_1x^2y+y^s,xy^3)
        \end{aligned}
    \end{equation}
    with $e_1 \in F$. That is, {$f \sim \widetilde{n}^2_{s,e_1}$.}
\end{proof}

\section{The classification of order 3}\label{class-3-allchar}
In this part we assume $f=(f_1,f_2) \in F[x,y]^2$ with {$\mathrm{ord}(f_1)=3,\ \mathrm{ord}(f_2) \geq 3.$} We also assume char$F=p>3$ in this part. We begin by classifying 3-jets.
\subsection{The classification of {3-jets}}

First choose a suitable coordinate system such that $j_3(f_1)=ax^3+bx^2y+cxy^2+dy^3$ with $a,b,c,d \in F$ and $a \neq 0.$ Then $j_3(f_1) \sim x^3+\frac{b}{a}x^2y+\frac{c}{a}xy^2+\frac{d}{a}y^3 \sim (x-e_1y)(x-e_2y)(x-e_3y) \sim l_1l_2l_3$ since $F$ is algebraically closed, where {$l_i,\ i=1,2,3,$} are linear forms in $R$.

I. $l_1=l_2=l_3$.  Let $\phi(x)=l_1,\ \phi(y)=y,$ then $\phi(j_3(f_1))=x^3$, i.e. $j_3(f_1) \sim x^3.$ We have $j_3(f_1,f_2) \sim (x^3, ay^3+bxy^2+cx^2y).$

\hspace*{1em}I.1 If $a=b=c=0$, $j_3(f) \sim (x^3,0).$

\hspace*{1em}I.2 If $a=b=0,\ c \neq 0,\ j_3(f) \sim (x^3, x^2y).$

\hspace*{1em}I.3 If $a=c=0,\ b \neq 0,\ j_3(f) \sim (x^3, xy^2).$

\hspace*{1em}I.4 {If $a=0,\ b,c \neq 0,\ j_3(f) \sim (x^3, x^2y+b'xy^2),$ where $b'=\frac{b}{c}$.} Let $\widetilde{x}=2b'x,\ \widetilde{y}=y$, then $j_3(f) \sim (\widetilde{x}^3, 2\widetilde{x}^2\widetilde{y}+\widetilde{x}\widetilde{y}^2) \sim (x^3, x^3+2x^2y+xy^2) \sim (x^3, x(x+y)^2).$ Let $\phi(x)=x,\ \phi(y)=y-x,$ then $j_3(f) \sim (x^3,xy^2).$ 

\hspace*{1em}I.5 If $a \neq 0,\ b^2 \neq 3ac$, we have 
\begin{equation}  {
    \begin{aligned}
        j_3(f) &\sim (x^3, y^3+b'xy^2+c'x^2y) \\
        &\sim (x^3,(y+\frac{b'}{3}x)^3+c'x^2y-\frac{b'^2}{3}x^2y-\frac{b'^3}{27}x^3)\\
        &\sim (x^3, (y+\frac{b'}{3}x)^3+(c'-\frac{b'^2}{3})x^2(y+\frac{b'}{3}x)),
    \end{aligned}}
\end{equation}
where $b'=\frac{b}{a},\ c'=\frac{c}{a}$. Using the automorphism $\phi(x)=x,\ \phi(y)=y-\frac{b'}{3}x$ we have $$j_3(f) \sim (x^3,y^3+(c'-\frac{b'^2}{3})x^2y).$$ Using the automorphism $\phi(x)=(c'-\frac{b'^2}{3})^{-\frac{1}{2}}x,\ \phi(y)=y$ since $c'-\frac{b'^2}{3} \neq 0$, we have $j_3(f) \sim (x^3,y^3+x^2y)$.

\hspace*{1em}I.6 If $a \neq 0,\ b^2=3ac$, as above, we have $j_3(f) \sim (x^3,y^3+(c'-\frac{b'^2}{3})x^2y) \sim (x^3,y^3).$

II. $l_1=l_2\neq l_3$.  Let $\phi(x)=l_1,\ \phi(y)=l_2,$ then $j_3(f_1) \sim x^2y.$ We have  {$j_3(f) \sim (x^2y, ax^3+bxy^2+cy^3)$.}

\hspace*{1em}II.1 If $a=b=c=0,\ j_3(f) \sim (x^2y,0).$

\hspace*{1em}II.2 If $a=b=0,\ c \neq 0,\ j_3(f) \sim (x^2y,y^3).$ Let $\phi(x)=y,\ \phi(y)=x,$ we get $j_3(f) \sim (x^3,xy^2).$

\hspace*{1em}II.3 If $b=c=0,\ a \neq 0,\ j_3(f) \sim (x^2y,x^3) \sim (x^3,x^2y).$

\hspace*{1em}II.4 If $a=c=0,\ b \neq 0,\ j_3(f) \sim (x^2y,xy^2).$

\hspace*{1em}II.5 If $c=0,\ a,b \neq 0,\ j_3(f) \sim (x^2y,ax^3+bxy^2)$. Let $\phi(x)=x,\ \phi(y)=\sqrt{\frac{a}{b}}y,$ we have $j_3(f) \sim (x^2y,x^3+xy^2).$

\hspace*{1em}II.6 If $c \neq 0,\ a \neq 0,$ we have
\begin{equation}
    \begin{aligned}
        j_3(f) &\sim (x^2y,y^3+a'xy^2+b'x^3) \\
        &\sim (x^2y,(y+\frac{a'}{3}x)^3+(b'-\frac{a'^3}{27})x^3),
    \end{aligned}
\end{equation}
where $a'=\frac{a}{c},\ b'=\frac{b}{c}$. Let $\phi(x)=x,\ \phi(y)=y+\frac{a'}{3}x,$ we have 
\begin{equation}
    \begin{aligned}
        j_3(f) &\sim (x^2(y-\frac{a'}{3}x),y^3+(b'-\frac{a'^3}{27})x^3) \\
        &\sim (x^2y-\frac{a'}{3}x^3,y^3+\frac{3(b'-\frac{a'^3}{27})}{a'}x^2y).
    \end{aligned}
\end{equation}
Let $\phi(x)=x,\ \phi(y)=-\frac{a'}{3}y,$ we have $j_3(f) \sim (x^3+x^2y,\ y^3+\lambda x^2y),$ where $\lambda=\frac{27b'-a'^3}{a'^3}=\frac{27bc^2-a^3}{a^3} \in F.$ If $\lambda=0$, it back to case I. Hence we assume $\lambda \neq 0$ here.

\hspace*{1em}II.7 If $c \neq 0,\ a=0,$ as above, $j_3(f) \sim (x^2(y-\frac{a'}{3}x),y^3+(b'-\frac{a'^3}{27})x^3) \sim (x^2y, y^3+b'x^3) \sim (x^2y, x^3+y^3).$

III. $l_1\neq l_2 \neq l_3$. By multiplying a unit, one can assume $l_3=\frac{1}{2}l_1+\frac{1}{2}l_2.$ Suppose $l_1=ux+vy,\ l_2=rx+sy$. Let $\phi(x)=\frac{u+r}{2}x+\frac{v+s}{2}y,\ \phi(y)=\frac{u-r}{2i}x+\frac{v-s}{2i}y,$ where $i^2=-1$. Then $\phi(l_1)=x+iy,\ \phi(l_2)=x-iy,\ \phi(l_3)=x$ and $j_3(f_1) \sim \phi(l_1l_2l_3)=x^3+xy^2.$ Hence $j_3(f) \sim (x^3+xy^2, axy^2+bx^2y+cy^3)$.

\hspace*{1em}III.1 If $a=b=c=0,\ j_3(f) \sim (x^3+xy^2,0).$

\hspace*{1em}III.2 If $b=c=0,\ a \neq 0,\ j_3(f) \sim (x^3,xy^2).$


\hspace*{1em}III.3 If $c=0,\ b \neq 0,\ j_3(f) \sim (x^3+xy^2, axy^2+bx^2y) \sim (x^3+xy^2, -ax^3+bx^2y) \sim (x^3+xy^2, x^2(y-\frac{a}{b}x))$. Let $\phi(x)=x,\ \phi(y)=y+\frac{a}{b}x,$ we have $j_3(f) \sim (x^2y, x^3+x(y+\frac{a}{b}x)^2)$, which goes back to II.

\hspace*{1em}III.4 If $c \neq 0,\ a=b=0,\ j_3(f) \sim (x^3+xy^2, y^3)$, which goes back to I.

\hspace*{1em}III.5 If $c \neq 0,$ one of $a,b \neq 0$, then write $j_3(f) \sim (x^3+xy^2,y^3+uxy^2+vx^2y)$, where $u,v \in F$ and one of $u,v \neq 0$. We have $j_3(f) \sim (x^3+xy^2,y^3+uxy^2+vx^2y+\alpha(x^3+xy^2))$. Choose $\alpha,r,s\in F$ such that 
\begin{equation}
    s \neq 0,\ u+\alpha=r+2s,\ v=2rs+s^2,\ \alpha=rs^2.
\end{equation}
These equations then reduced to be
\begin{equation}\label{III-3}
    s^4+(3-v)s^2-2us-v=0,\ 2r=\frac{v}{s}-s,\ \alpha=rs^2, 
\end{equation}
hence such $\alpha,r,s$ {\color{black}exist}. Then $j_3(f) \sim (x^3+xy^2,(y+rx)(y+sx)^2)$. Using the automorphism $\phi(x)=x+\frac{y}{s},\phi(y)=y+rx$, {\color{black}then it is deduced to case I or II.}

Hence we get the result:

\begin{proposition}
    Let $f \in F[x,y]^2$ with $\mathrm{ord}(f)=3$ be a unimodal complete intersection singularity, then $j_3(f)$ is equivalent to one of the following:
\begin{equation}  {
    \begin{split}
        &(x^3,0),\ (x^3,x^2y),\ (x^3,xy^2),\ (x^3,y^3+x^2y),\ (x^3,y^3),\\
        & (x^2y,0),\ (x^2y,xy^2),\ (x^2y,x^3+xy^2),\ (x^3+x^2y,y^3+ \lambda x^2y)(\lambda \neq 0),\\
        &(x^2y,x^3+y^3),\ (x^3+xy^2,0).
    \end{split}}
\end{equation}

\end{proposition}

\subsection{The classification of unimodal}

We have the following classification of unimodal ICIS of order 3:
\begin{proposition}\label{ord3-class}
    {A unimodal ICIS of order 3 in any field with characteristic not equal to $2,3$ must have the form in Table \ref{table-ord3}.}

\renewcommand\arraystretch{1.5}
\begin{longtable}{|c|c|c|}
\caption{}\label{table-ord3}
\\
\hline
Symbol&Form& condition  \\
\hline
$H$&$(x^3+x^2y,y^3+\lambda x^2y)$&$\lambda \neq 0$ \\
\hline
$I$&$(x^3,y^3+x^2y)$& \\
\hline
$J$&$(x^3,y^3)$ & \\
\hline
$K_r$&$(x^3+xy^2,x^2y+ y^r)$&$r \geq 4$ \\
\hline
$L_{r,s}$&$(x^2y+y^r,xy^2+x^s)$&$r,s \geq 4$\\
\hline
$M_r$&$(x^3+y^r,xy^2)$ & {$r \geq 4$}\\
\hline
$N_{\lambda}$&$(x^3+\lambda xy^3,x^2y+y^4)$&$ \lambda \notin \{1,12\}$\\
\hline
$P_{r,\infty}$ & $(x^3+xy^3+xy^r,x^2y+y^4)$ & $r \geq 4$\\
\hline
$P_{\infty,s}$ & $(x^3+xy^3+y^s,x^2y+y^4)$ & $s \geq 5$ \\
\hline
$P_{r,s,\lambda}$ & $(x^3+xy^3+xy^r+\lambda y^s,x^2y+y^4)$ & $r \geq 4,s \geq 5,$ \\
&& $\lambda \in F$\\
\hline
$\widetilde{P}_{r,s} $&$(x^3+xy^3+uxy^r+vy^s,x^2y+y^4),\ where$& {$r \geq 4, \ s \geq 5,$}\\
&$u=u_0+u_1y+\dots,\ {v=v_0+v_1y+\dots,}$& $p \mid 2r-2s+3$\\
\hline
$R_t$ &$(x^3+xy^3,x^2y+y^t)$&$ t \geq 5$\\
\hline
$X_{\lambda}$&$(x^3+y^4,x^2y+\lambda y^4)$&$\lambda \in \{0,1\}$\\
\hline
$Y_{\lambda}$&$(x^3+y^5,x^2y+\lambda y^5)$&$\lambda \in \{0,1\}$\\
\hline
$Z_{\lambda}$&$(x^3+12xy^3+\lambda y^5,x^2y+y^4)$&$\lambda \in \{0,1\}$\\
\hline
\end{longtable}
\end{proposition}

We will prove {Proposition \ref{ord3-class}} step by step:

\begin{proposition}\label{compute-ex2}
    If $j_3(f)$ is contact equivalent to one of {\color{black}the following forms}
    \begin{equation}\label{finite-determined-jet}
        (x^3,y^3+x^2y),\ (x^3,y^3),\ (x^2y,x^3+y^3),\ (x^3+x^2y,y^3+ \lambda x^2y)(\lambda \neq 0),
    \end{equation}
    then $j_3(f)$ is 3-determined. In particular, $f$ is contact equivalent to   $j_3(f)$ of the above forms.
\end{proposition}
\begin{proof}
    After some computation, one can show $\mathfrak{m}^4 \subset \mathfrak{m}\cdot \widetilde{T}_f(\mathcal{K}f)$ when $j_3(f)$ has forms in  {(\ref{ord3-class})}. Hence by Theorem \ref{finite-determinacy}, $j_3(f)$ is 3-determined. Here we compute the case when $j_3(f) \sim (x^3,y^3+x^2y)$ as an example.

    Let
$$e_1=(x^3,0),\ e_2=(0,x^3),\ e_3=(y^3+x^2y,0),$$ 
       $$ e_4=(0,y^3+x^2y),\ e_5=(3x^2,2xy),\ e_6=(0,3y^2+x^2).$$
    Then $\widetilde{T}_f(\mathcal{K}f)$ is generated by $e_1,\ e_2,\ e_3,\ e_4,\ xe_5,\ ye_5,\ xe_6,\ ye_6.$ And we have:\\
    $(0,x^4)=xe_2 \in \mathfrak{m}\widetilde{T}_f(\mathcal{K}f).$\\
    $(0,x^3y)=ye_2 \in \mathfrak{m}\widetilde{T}_f(\mathcal{K}f).$\\
    $(0,x^2y^2)=\frac{1}{3}(x^2e_6-xe_2) \in \mathfrak{m}\widetilde{T}_f(\mathcal{K}f).$\\
    $(0,xy^3)=\frac{1}{3}(xye_6-(0,x^3y)) \in \mathfrak{m}\widetilde{T}_f(\mathcal{K}f).$\\
    $(0,y^4)=ye_4-(0,x^2y^2) \in \mathfrak{m}\widetilde{T}_f(\mathcal{K}f).$\\
    $(x^4,0)=xe_1 \in \mathfrak{m}\widetilde{T}_f(\mathcal{K}f).$\\
    $(x^3y,0)=ye_1 \in \mathfrak{m}\widetilde{T}_f(\mathcal{K}f).$\\
    $(x^2y^2,0)=\frac{1}{3}(y^2e_5-2(0,xy^3)) \in \mathfrak{m}\widetilde{T}_f(\mathcal{K}f).$\\
    $(xy^3,0)=xe_3-(x^3y,0) \in \mathfrak{m}\widetilde{T}_f(\mathcal{K}f).$\\
    $(y^4,0)=ye_3-(x^2y^2,0) \in \mathfrak{m}\widetilde{T}_f(\mathcal{K}f).$\\
    This implies $\mathfrak{m}^4 \cdot R^2 \subset \mathfrak{m}\cdot \widetilde{T}_f(\mathcal{K}f)$, as claimed.
\end{proof}

\begin{proposition}
    If $j_3(f) \sim (x^3,0),(x^2y,0),(x^3+xy^2,0)$, then $\mathcal{K}-mod(f)$ is at least 2. 
\end{proposition}
\begin{proof}
    We just prove the case when $j_3(f) \sim (x^2y+xy^2,0)$. The others are similar.

    In this case a complete transversal in $J_4$ is spanned by $(0,x^4),(0,x^2y^2),(0,y^4)$, hence $j_4(f) \sim (x^3+xy^2,ax^4+bx^2y^2+cy^4)$ by Theorem \ref{complete-cor}. Let $g=(x^3+xy^2,ax^4+bx^2y^2+cy^4)$. After computation, we have $(0,xy^3),(0,y^4) \notin P_{3,4}$ for almost all $a,b,c$. Hence, the codimension of $\widetilde{T}_g(\mathcal{K}_4 g)$ in $P_{3,4}$ is 2, which is, the modality of $f \geq 2$ by Proposition \ref{cod-modality}.
\end{proof}

\begin{proposition}\label{x3-xy2}
    If $j_3(f) \sim (x^3,xy^2)$, then $f \sim (x^3+y^r,xy^2),\ r \geq 4.$
\end{proposition}
\begin{proof}
    We can compute the complete transversal as follows: By Proposition \ref{tangent-image}, $\widetilde{T}_f(\mathcal{K}f)$ is generated by $$e_1=(x^3,0),\ e_2=(0,x^3),\ e_3=(xy^2,0),\ e_4=(0,xy^2),\ xe_5,\ ye_5,\ xe_6,\ ye_6,$$ where $e_5=(3x^2,y^2),\ e_6=(0,2xy)$. A similar computation shows that $$(x^4,0),\ (x^3y,0),\ (x^2y^2,0),\ (xy^3,0),\ (0,x^4),\ (0,x^3y),\ (0,x^2y^2),\ (0,xy^3),\ (0,y^4) $$ lie in $\widetilde{T}_f(\mathcal{K}f)$ while $(y^4,0) \notin \widetilde{T}_f(\mathcal{K}f)$. In fact, one can easily show that $(y^l,0) \notin \widetilde{T}_f(\mathcal{K}f)$ for any $l>3$. Hence, a complete transversal is spanned by $\{ (y^l,0) \mid l \geq 4\}$.

    By {Theorem \ref{complete-cor}}, we have 
    \begin{equation} {
        \begin{aligned}
            f &\sim (x^3+\sum_{l\geq 4}b_l y^l,xy^2) \\
            &\sim (x^3+e(y)y^r,xy^2),\ r \geq 4,\ e(y)\ \mathrm{is\ a\ unit.} \\
            &\sim (x^3+y^r,xy^2).
        \end{aligned} }   
    \end{equation}
In the last line, we take the automorphism $\phi(x)=e(y)^{\frac{1}{3}}x$ and $\phi(y)=y.$
\end{proof}

\begin{proposition}\label{rs-trick}
    If $j_3(f) \sim (x^2y,xy^2),$ then $f \sim (x^2y+y^r,xy^2+x^s),\ r,s \geq 4.$
\end{proposition}
\begin{proof}
    Use a computation like Proposition \ref{x3-xy2}, a complete transversal is given by $\{(y^l,0),(0,x^l) \mid l\geq 4\}.$ Hence by {Theorem} \ref{complete-cor}, we have

    \begin{equation}
        \begin{aligned}
            f &\sim (x^2y+\sum_{l\geq 4}a_l y^l,xy^2+\sum_{l\geq 4}b_l x^l) \\
            &\sim (x^2y+a(y)y^r,xy^2+b(x)x^s),
        \end{aligned}
    \end{equation}
    where $a(y),b(x)$ are units in $F[[x,y]]$ and $r,s \geq 4$.

    Using the automorphism $\phi(x)=a(y)^{\frac{1}{2}}x$ and $\phi(y)=y$, we have 
\begin{equation}
    \begin{aligned}
        f &\sim (x^2y+y^r,a(y)^{\frac{1}{2}}xy^2+a(y)^{\frac{s}{2}}b(a(y)^{\frac{1}{2}}x)x^s)\\ & \sim (x^2y+y^r,xy^2+e(x,y)x^s),
    \end{aligned}
\end{equation}
     where $e(x,y)$ is a unit. We write $e(x,y)=\sum_{i \geq 0}e_i(x)y^i,$ then
     \begin{equation}
         \begin{aligned}
             f &\sim (x^2y+y^r,xy^2+(\sum_{i \geq 0}e_i(x)y^i)x^s) \\ 
             & \sim (x^2y+y^r,xy^2+(\sum_{i \geq 0}e_i(x)y^i)x^s-x^{s-2}(x^2y+y^r)(\sum_{k \geq 1} e_k(x)y^{k-1}))\\
             &\sim (x^2y+y^r,xy^2+e_0(x)x^s-\sum_{k \geq 1}x^{s-2}{e_{k}(x)}y^{r+k-1}).
         \end{aligned}
     \end{equation}
     {Note that the order of $x$ changes from $s$ to $s-2$, and the order of $y$ changes from $0$ to $r$. Repeat the operation, we get $f \sim (x^2y+y^r,xy^2+e_0(x)x^s+d(x,y))$, where $d(x,y)=y^{\frac{rs}{2}}e'(x,y)$ or $xy^{\frac{r(s-1)}{2}}e'(x,y)$ depending on $s$ is even or odd, and $e'(x,y)=\sum_{k \geq 1}e_k(x)y^{k-1}$, which is a unit in $F[[x,y]]$.} The order of $d$ is $\geq \frac{r(s-1)}{2}+1$.

    Using the automorphism $\phi(x)=x$ and $\phi(y)=e_0(x)^{\frac{1}{2}}y$. Then  
    \begin{equation}
        \begin{aligned}
            f &\sim (e_0(x)^{\frac{1}{2}}x^2y+e_0(x)^{\frac{r}{2}}y^r,xy^2+x^s+d_1(x,y))\\
            &\sim (x^2y+e_0(x)^{\frac{r-1}{2}}y^r,xy^2+x^s+d_1(x,y)),
        \end{aligned}
    \end{equation}
    {where $d_1(x,y)=d(x,e_0(x)^{\frac{1}{2}}y)$.}
    Write $e_0(x)^{\frac{r-1}{2}}=\sum_{i \geq 0}u_i x^i,$ then 
    \begin{equation}  {
        \begin{aligned}
            f &\sim (x^2y+(\sum_{i \geq 0}u_i x^i)y^r,xy^2+x^s+d_1(x,y))\\
            &\sim (x^2y+(\sum_{i \geq 0}u_i x^i)y^r-y^{r-2}(xy^2+x^s+d_1(x,y))(\sum_{k \geq 1}u_k x^{k-1}),xy^2+x^s+d_1(x,y))\\
            &\sim (x^2y+u_0y^r-\sum_{k \geq 1}u_k y^{r-2}x^{s+k-1}-\sum_{k \geq 1}u_k y^{r-2}d_1(x,y),xy^2+x^s+d_1(x,y)).
        \end{aligned}}
    \end{equation}
    Repeat the operation, we get  {$f \sim (x^2y+u_0y^r+d_2(x,y),xy^2+x^s+d_1(x,y)),$} where the order of {$d_2(x,y)$} is $ \geq \frac{s(r-1)}{2}+1.$ 
    
    Taking the automorphism $\phi(x)=\alpha x,\ \phi(y)=\beta y,$ where $\alpha,\beta \in F$ satisfy $\alpha^2\beta=u_0 \beta^r,\ \alpha\beta^2=\beta^s$, we have {$$f \sim (x^2y+y^r+\widetilde{d_2}(x,y),xy^2+x^s+\widetilde{d_1}(x,y)).$$} Specifically, {$\widetilde{d_1}(x,y)$} has order $\frac{rs}{2}$ if $s$ is even and $\frac{r(s-1)}{2}+1$ if $s$ is odd. {$\widetilde{d_2}(x,y)$} has order $\frac{rs}{2}$ if $r$ is even and $\frac{(r-1)s}{2}+1$ if $r$ is odd.

    Now we exchange the position of $x,y$ so that $r \geq s$. Let $g=j_r(f)=(x^2y+y^r,xy^2+x^s)$. Using the similar computation {as the proof of \textbf{(3.b.ii.2.3.3.1)} in Proposition \ref{ord2-class-mod1}}, we can show $\mathfrak{m}^{r+1}R^2 \subset \mathfrak{m} \widetilde{T}_g(\mathcal{K}g)$. This means $g$ is {$(2r-3)$}-determined by Theorem \ref{finite-determinacy}. Since $s \geq 4$, {$min\{\mathrm{ord}\widetilde{d_1}(x,y),\mathrm{ord}\widetilde{d_2}(x,y)\}>2r-3$}. Therefore $f \sim j_{2r-3}(f)=(x^2y+y^r,xy^2+x^s)$.

\end{proof}

\begin{proposition}\label{x2y-x3+xy2}
    If $j_3(f) \sim (x^2y,x^3+xy^2),$ then $f \sim (x^2y+y^r,x^3+xy^2),\ r \geq 4.$
\end{proposition}
\begin{proof}
    A complete transversal is given by $\{(y^l,0) \mid l\geq 4\}$, hence we have
    $$f \sim (x^2y+\sum_{l\geq 4}a_l y^l,x^3+xy^2) \sim (x^2y+e(y)y^r,x^3+xy^2),$$ where $r \geq 4.$

    Using the automorphism $\phi(x)=e(y)^{\frac{1}{2}}x$ and $\phi(y)=y,$ we have $$f \sim (x^2y+y^r,e(y)x^3+xy^2).$$ Write $e(y)=\sum_{i \geq 0} e_iy^i$, then 
    \begin{equation}
        \begin{aligned}
            f &\sim (x^2y+y^r,(\sum_{i \geq 0} e_iy^i)x^3+xy^2-x(x^2y+y^r)(\sum_{k \geq 1}e_ky^{k-1})) \\
            &\sim (x^2y+y^r,e_0x^3+xy^2-\sum_{k \geq 1}e_kxy^{r+k-1}) \\
            &\sim (x^2y+y^r,e_0x^3+xy^2(1-\sum_{k \geq 1}e_kxy^{r+k-3})).
        \end{aligned}
    \end{equation}
    Denote $u_0(y)=e(y),\ v_1(y)=1-\sum_{k \geq 1}e_kxy^{r+k-3}.$ Then $v_1(y) \in 1+\mathfrak{m}^{r-2}$ and 
    \begin{equation}
        \begin{aligned}
            f &\sim (x^2y+y^r,e_0x^3+u_1(y)xy^2) \\
            &\sim (x^2y+y^r,e_0v_1(y)^{-1}x^3+xy^2).
        \end{aligned}
    \end{equation}
    Since $v_1(y) \in 1+\mathfrak{m}^{r-2},$ we have $v_1(y)^{-1} \in 1+\mathfrak{m}^{r-2}$. Denote $u_1(y)=e_0v_1(y)^{-1},$ then $u_1(y) \in e_0+\mathfrak{m}^{r-2}$ and $f \sim (x^2y+y^r,u_1(y)x^3+xy^2)$. Write {$$u_1(y)=e_0+\sum_{k \geq r-2}v_k y^k,$$}we have 
    \begin{equation}
        \begin{aligned}
            f &\sim (x^2y+y^r,(e_0+\sum_{k \geq r-2}v_k y^k)x^3+xy^2)\\
            &\sim (x^2y+y^r,(e_0+\sum_{k \geq r-2}v_k y^k)x^3+xy^2-x(x^2y+y^r)(\sum_{k \geq r-2}v_k y^{k-1})) \\
            &\sim (x^2y+y^r,e_0x^3+xy^2-\sum_{k \geq r-2}v_kxy^{r+k-1}) \\
            &\sim (x^2y+y^r,e_0x^3+xy^2(1-\sum_{k \geq r-2}v_ky^{r+k-3})).
        \end{aligned}
    \end{equation}
    Denote $v_2(y)=1-\sum_{k \geq r-2}v_ky^{r+k-3},$ then $v_2(y) \in 1+\mathfrak{m}^{2r-5}.$ Denote $u_2(y)=e_0v_2(y)^{-1} \in e_0+\mathfrak{m}^{2r-5},$ we have $f \sim (x^2y+y^r,u_2(y)x^3+xy^2).$ 
    
    Repeat the operation, we can get a sequence of units $u_1(y),u_2(y), \dots, u_n(y), \dots$ with $u_n(y) \in e_0+\mathfrak{m}^{n(r-3)+1}$ and $f \sim (x^2y+y^r,u_n(y)x^3+xy^2)$. 
    
    Since $r \geq 4,$ the orders of $u_n(y)-e_0$ are strictly increasing. Note that $(x^2y+y^r,e_0x^3+xy^2)$ has finite Tjurina number (and hence is finite determined), we have $f \sim (x^2y+y^r,e_0x^3+xy^2) \sim (x^2y+y^r,x^3+xy^2)$ by using Remark \ref{finite-deter-tau} and applying the automorphism $\phi(x)=\alpha x,\ \phi(y)=\beta y$ where $\alpha^2 \beta=\beta^r$ and $e_0\alpha^3=\alpha \beta^2.$
\end{proof}

The situation becomes complicated when $j_3(f) \sim (x^3,x^2y).$ A complete transversal is given by $\{(xy^{l-1},0),(y^l,0),(0,y^l) \mid l \geq 4\}$. Then 
\begin{equation}\label{big-f}
\begin{aligned}
    f &\sim (x^3+\sum_{i \geq 3}a_i xy^{i}+\sum_{j \geq 4}b_jy^j,x^2y+\sum_{k \geq 4}c_ky^k)\\
    & \sim (x^3+a(y)xy^r+b(y)y^s,x^2y+c(y)y^t),
\end{aligned}
\end{equation}
where $r \geq 3,\ s \geq 4,\ t \geq 4,$ $a(y),b(y),c(y)$ are units or 0.

First we have the following criterion:

\begin{proposition}
    The modality of $f \sim (x^3+a(y)xy^r+b(y)y^s,x^2y+c(y)y^t) $ is at least $2$ if $r \geq 4,\ s \geq 6$ and $t \geq 5.$
\end{proposition}
\begin{proof}
    If $r \geq 4,\ s \geq 6$ and $t \geq 5,$ set $g=j_4(f)=(x^3,x^2y)$, then all jets in $J_5'(g)$ are equivalent to $g_{ac}:= (x^3+axy^4,x^2y+cy^5),$ where $J_5'(g)$ is formed by jets in $J_5(g)$ with $r \geq 4,\ s \geq 6$ and $t \geq 5.$ 
    
    {Analogously to} Proposition \ref{ord2-modality1}, we can show $(y^5,0),(0,y^5) \notin \widetilde{T}_{g_{ac}}(\mathcal{K}_5 g_{ac})$. Hence $\mathrm{cod}(g_{ac}) \geq 2$. By Proposition \ref{cod-modality} and Proposition \ref{modality}, $\mathcal{K}-mod(f) \geq \mathcal{K}_5-mod(f) \geq 2.$
\end{proof}

\begin{remark}
    If $a(y)\ (resp.\ b(y),c(y))=0$, we regard $r\ (resp.\ s,t)=\infty$ here.
\end{remark}

Now we assume $r \leq 4$ or $s \leq 6$ or $t \leq 5$ and $a(y),b(y),c(y)$ are units. We can simplify {(\ref{big-f})} by the trick using in \cite{icissimple}.

First we recall the Implicit Function Theorem:

\begin{lemma}\label{implicit-function}
    (cf. \cite{singular-introduction} Theorem 6.2.17, Implicit Function Theorem) Let $\mathcal{K}$ be a field and $F \in \mathcal{K}[[x_1,\dots,x_n,y]]$ such that
    \begin{equation}
        F(x_1,\dots,x_n,0) \in \langle x_1,\dots,x_n \rangle,\ \frac{\partial F}{\partial y}(x_1,\dots,x_n,0) \notin \langle x_1,\dots,x_n \rangle,
    \end{equation}
    then there exists a unique $y(x_1,\dots,x_n) \in \langle x_1,\dots,x_n \rangle \mathcal{K}[[x_1,\dots,x_n]]$ such that $$F(x_1,\dots,x_n,y(x_1,\dots,x_n))=0.$$
\end{lemma}

Back to $f$ in {(\ref{big-f})}, if 
\begin{equation}\label{condition-1}
    3r-2s \neq 0 \mathrm{\ and\ } p \nmid 3r-2s,
\end{equation}
{applying} $\phi(x)=a(y)^{\frac{1}{2}}x$ and $\phi(y)=y,$ we get $$f \sim (x^3+xy^r+y^s\widetilde{b}(y),x^2y+y^t c_1(y)),$$ where $\widetilde{b}(y)=\frac{b(y)}{a(y)},\ c_1(y)=\frac{c(y)}{a(y)}$.

Write $\widetilde{b}(y)=\sum_{i \geq 0}b_i y^i$. Consider the function
$$F(z)=z^{2s-3r}\sum_{i \geq 0}b_iy^iz^{2i}-b_0.$$
We have $F(1)\in \langle y \rangle F[[y]]$, and {$$\ F'(1)=(2s-3r)\sum_{i \geq 0}b_iy^i-2\sum_{i \geq 1}i b_iy^i $$} is a unit since $3r-2s \neq 0 \mathrm{\ and\ } p \nmid 3r-2s$. Apply Lemma \ref{implicit-function} to function $G(z)=F(z+1)$, there exists a $\widetilde{z}(y)$ such that $G(\widetilde{z}(y))=0.$ Let $z(y)=\widetilde{z}(y)+1$ and then $z(y)$ is a unit and $F(z(y))=0$.

Using the automorphism $\phi(x)=z(y)^r x$ and $\phi(y)=z(y)^2 x$, we have $f \sim (x^3+xy^r+b_0y^s,x^2y+y^t c_2(y))$. Then apply $\xi (x)=\alpha x,\ \xi(y)=\beta(y)$ with $\alpha,\ \beta \in F$ satisfying $\alpha^3=\alpha\beta^r,\ \alpha\beta^r=b_0\beta^s$ (such $\alpha,\ \beta$ exists since $3r-2s \neq 0$), we have {
\begin{equation}\label{form-1}
    f \sim (x^3+xy^r+y^s,x^2y+y^t \widetilde{c}(y)),
\end{equation}
where $\widetilde{c}(y)$ is the image of $c_2(y)$ under the automorphism $\xi$.}

Similarly, if
\begin{equation}\label{condition-2}
    2s-3t+3 \neq 0 \mathrm{\ and\ } p \nmid 2s-3t+3,
\end{equation}
we have {
\begin{equation}\label{form-2}
    f \sim (x^3+\widetilde{a}(y)xy^r+y^s,x^2y+y^t),
\end{equation}
where $\widetilde{a}(y)$ is the image of $a(y)$ under the similar automorphism.}
 
 If
 \begin{equation}\label{condition-3}
     r+1-t \neq 0 \mathrm{\ and\ } p \nmid r+1-t,
 \end{equation}
 we have
\begin{equation}\label{form-3}
    f \sim (x^3+xy^r+\widetilde{b}(y)y^s,x^2y+y^t).
\end{equation}

If all {(\ref{condition-1}), (\ref{condition-2}), (\ref{condition-3})} fail, then $r,s,t$ satisfy {
\begin{equation}\label{all-condition} 
    \begin{aligned}
        3r-2s &=ap,\\
        2s-3t+3&=bp,\\
        r+1-t&=cp.
    \end{aligned}
\end{equation}}
But the minimal solution of {(\ref{all-condition})} with $r\geq 3,\ s \geq 4,\ t \geq 4$ and $p \geq 5$ is exactly $r=4,\ s=6,\ t=5$, in which case $f$ is not unimodal.

\begin{remark}\label{alpha-beta}
    (i) We call the technique we use here as $\alpha,\beta$-trick, since we can easily (but not strictly) apply $\phi(x)=\alpha x,\phi(y)=\beta x$ in {(\ref{big-f})} and get the result. For example, when {(\ref{condition-1})} holds, choose $\alpha,\beta$ as the solution of {
    \begin{equation}
        \begin{cases}
\alpha^3=\alpha \beta^r a(y),\\
\alpha \beta^r a(y)=\beta^s b(y),
\end{cases}
    \end{equation}}
    then apply $\phi$ on {(\ref{big-f})}, we can get {(\ref{form-1})}.
    The trick was shown in \cite{complete-unimodal-plane}. The implicit function theorem provides a complete proof with a same result. But it is useless in some special characteristic, e.g. $p \mid 3r-2s$ here.

    (ii) We call the technique used in Proposition \ref{rs-trick} and \ref{x2y-x3+xy2} as $r,s$-trick. It has no restriction on characteristic, although the process is a little tedious. Later we will use the $r,s$-trick again but omit the process.
\end{remark}

\begin{proposition}\label{s4}
    If $s=4$, then {$f \sim (x^3+y^4,x^2y+\lambda y^4) \sim X_{\lambda}$ in Table \ref{table-ord3}}, where $\lambda=0 \mathrm{\ or\ }1$.
\end{proposition} 
\begin{proof}
    Let $h=(x^3+y^4,x^2y)$ be the weighted 0-jet of $f$ with respect to $(a,d)$, where $a=(4,3),d=(12,11).$ 
    
    We choose $T=\mathrm{span}\langle (0,y^4) \rangle$ as a complete transversal. Then we have $F_{a,d}^1 \subset T+\widetilde{T}_h(F_{a,d}^1\mathcal{K}h).$ By Proposition \ref{weight-ct}, $f \sim (x^3+y^4,x^2y+ay^4)$ for some $a \in F$. After an obvious scaling, $f \sim (x^3+y^4,x^2y+y^4)$ or $(x^3+y^4,x^2y)$.

\end{proof}

\begin{proposition}\label{s5r3}
    If $s=5$ and $r=3$, then {$f \sim N_{\lambda},R_t,P_{r,s},Z_{\lambda}$ in Table \ref{table-ord3}.} 
\end{proposition}
\begin{proof}
    In this case $p \nmid 3r-2s$ holds. By {(\ref{form-1})}, we have $j_4(f) \sim (x^3+xy^3,x^2y+\lambda_0 y^4)$, which is weighted homogeneous of degree 0 with respect to $(a;d)$, where $a=(3,2),d=(9,8)$. Write $g=j_4(f)$. Next we find $T \subset F_{a,d}^1R^2\backslash F_{a,d}^rR^2$ such that {$F_{a,d}^1R^2 \subset T+\widetilde{T}_g(F^1\mathcal{K}g)$}. Then by Proposition \ref{weight-ct}, $f \sim g+t$ with $t \in T$.

    After some computation (can be easily done by hand), we know that for $\lambda_0 \neq 0, 1, \frac{1}{12}$, $F_{a,d}^1R^2 \subset \widetilde{T}_g(F^1\mathcal{K}g)$ (note that $charF \neq 2,3$ hence $\frac{1}{12}$ is well-defined). In that case $f \sim g = (x^3+xy^3,x^2y+\lambda_0 y^4) \sim (x^3+\lambda xy^3,x^2y+ y^4)$, where $\lambda=\frac{1}{\lambda_0}$.

    If $\lambda_0=0$, then $F_{a,d}^1R^2 \subset T+\widetilde{T}_g(\mathcal{K}\cdot g)$ where $T$ is spanned by $\{(0,y^5),(0,y^6),\dots\}$. Hence $f \sim (x^3+xy^3,x^2y+e y^t)$ with $t \geq 5$ or $f \sim (x^3+xy^3,x^2y)$ (which is not an ICIS). Using $r,s$-trick we can show $f \sim (x^3+xy^3,x^2y+y^t).$
    
    If $\lambda_0=1, \frac{1}{12}$, we write $g \sim (x^3+\lambda xy^3,x^2y+ y^4)$ with $\lambda=1,12$. The process will be shown later in Proposition \ref{s6t4}. 

\end{proof}

\begin{proposition}\label{s5r4}
    If $s=5$ and $r \geq 4$, then {$f \sim N_0$ or $Y_{\lambda}$ in Table \ref{table-ord3},} where $\lambda \in \{0,1\}$.
\end{proposition}
\begin{proof}
    We have $f \sim (x^3+a(y)xy^r+b(y)y^5,x^2y+c(y)y^t)$. 

    If $t=4$, let $g=j_4(f)=(x^3,x^2y+cy^4) \sim (x^3,x^2y+y^4)$ be the 4-jet of $f$. Then $g$ is weighted homogeneous of degree 0 with respect to $(a;d)$, where $a=(3,2),d=(9,8)$ as in Proposition \ref{s5r3}. After the same computation, we have {$F_{a,d}^1R^2 \subset \widetilde{T}_g(F^1\mathcal{K}g)$}. That is, $f \sim (x^3,x^2y+y^4)$.

    If $t \geq 5$, let $h=(x^3+by^5,x^2y) \sim (x^3+y^5,x^2y)$ be the weighted 0-jet with respect to $(a,d)$, where $a=(5,3),d=(15,13)$. Computation shows {$F_{a,d}^1R^2 \subset T+\widetilde{T}_g(F^1\mathcal{K}g)$} for $T=span\langle (0,y^5)\rangle$ (whether $p=5$ or not). Hence $f \sim (x^3+y^5,x^2y+ay^5)$. An obvious scaling shows $f \sim (x^3+y^5,x^2y+\lambda y^5)$ for $\lambda \in \{0,1\}$.
\end{proof}

\begin{proposition}\label{s6t4}
    If $s \geq 6$ and $t=4$, then {$f \sim N_\lambda,P_{r,\infty},P_{\infty,s},P_{r,s,\lambda},\widetilde{P}_{r,s}$  or $Z_\lambda$ in Table \ref{table-ord3}.}
\end{proposition}
\begin{proof}
    Consider the $4$-jet $j_4(f)=(x^3+axy^3,x^2y+cy^4) \sim (x^3+\lambda xy^3,x^2y+ y^4)$, where $\lambda =\frac{a}{c}$. If $\lambda \neq 1,12,$ a computation similar with Proposition \ref{s5r3} shows $j_4(f) \sim (x^3+\lambda xy^3,x^2y+ y^4)$ and furthermore $f \sim  (x^3+\lambda xy^3,x^2y+ y^4)$.

    If $\lambda=12$, we write $g=j_4(f)=(x^3+ 12xy^3,x^2y+ y^4)$, which is weighted homogeneous of degree 0 with respect to $(a,d)$, where $a=(3,2),d=(9,8)$. Computation shows that a weighted complete transversal is given by $(y^5,0)$, i.e. {$$F_{a,d}^1R^2 \subset span\langle (y^5,0) \rangle+\widetilde{T}_g(F^1\mathcal{K}g).$$} Then by Proposition \ref{weight-ct}, $f \sim (x^3+ 12xy^3,x^2y+ y^4)$ or $(x^3+ 12xy^3+y^5,x^2y+ y^4)$.

    If $\lambda_0=1$, we write $g \sim (x^3+ xy^3,x^2y+ y^4)$. Then {$F_{a,d}^1R^2 \subset T+\widetilde{T}_g(F^1\mathcal{K}g)$} with {$T=span\langle (y^l,0),(xy^j,0)|l \geq 5,j \geq 4 \rangle$.} Hence $f \sim (x^3+ xy^3+u(y)xy^r+v(y)y^s, x^2y+y^4)$ with $r \geq 5,s \geq 4$.

    If $v(y)=0,$ write $u(y)=u_0+u_1y+\dots.$ Through the process
    \begin{equation}\label{reduce-u-v-process}
        \begin{aligned}
            f &\sim (x^3+xy^3+u(y)xy^r,x^2y+y^4) \\
            &\sim ((x^3+xy^3+u(y)xy^r)(1-\frac{u_1}{u_0}y),x^2y+y^4) \\
            &\sim (x^3+xy^3+(u_0+(u_2-\frac{u_1^2}{u_0})y^2+\dots)xy^r-\frac{u_1}{u_0}x(x^2y+y^4),x^2y+y^4)\\
            &\sim (x^3+xy^3+(u_0+(u_2-\frac{u_1^2}{u_0})y^2+\dots)xy^r,x^2y+y^4),
        \end{aligned}
    \end{equation}
    we reduce $u_1$ to $0$. Repeat the process, we can reduce $u(y)$ to $u_0 \in F$ and finally to $1$. Then we have $f \sim (x^3+xy^3+xy^r,x^2y+y^4) \sim P_{r,\infty}$.

    If $u(y)=0$, similarly we can get $f \sim (x^3+xy^3+y^s,x^2y+y^r) \sim P_{\infty,s}$.

    If {$u(y),v(y) \neq 0$}, apply $\phi(x)=\alpha(y)^{-3},\phi(y)=\alpha(y)^{-2}$, we have $ f \sim (x^3+xy^3+u(y)\alpha(y)^{2r-6}xy^r+v(y)\alpha(y)^{2s-9}y^s,x^2y+y^s)$. Using $\alpha,\beta$-trick, if $p \nmid 2r-2s-3$, then there exists $\alpha(y)$ such that $\alpha(y){2r-6}u(y)=\alpha(y)^{2s-9}v(y)$. By the same process to {(\ref{reduce-u-v-process})}, we have $f \sim (x^3+xy^3+u_0xy^r+u_0y^s,x^2y+y^4) \sim (x^3+xy^3+xy^r+\lambda y^s,x^2y+y^4)\sim P_{r,s,\lambda}$.

    If $p \mid 2r-2s-3$, we get a family $\widetilde{P}_{r,s}$.
\end{proof}

\begin{proposition}\label{s-geq-6}
    If $r=3,s \geq 6, t \geq 5$, then {$f \sim R_t$ in Table \ref{table-ord3}.}
\end{proposition}

\begin{proof}
    In this case, let $g=j_4(f)=(x^3+xy^3,x^2y)$. Then an ordinary complete transversal is given by $T=span \langle (0,y^5),(0,y^6), \dots \rangle$. Hence $f \sim (x^3+xy^3,x^2y+e(y)y^t)$. Using $r,s$-trick we get $f \sim (x^3+xy^3,x^2y+y^t),\ t \geq 5.$
\end{proof}

The above propositions finish the proof of {Proposition \ref{ord3-class}.}

\section{The classification of order 3 when $char F=2$}

The process of classification in the field of characteristic 2 is quite similar to that of other characteristics. We first classify 3-jets and then classify all germs.

\subsection{The classification of {3-jets}}
Same as in Section \ref{class-3-allchar}, for $f=(f_1,f_2)$ with $\mathrm{ord}(f_1)=3$, we have $j_3(f_1) \sim (x-e_1y)(x-e_2y)(x-e_3y) \sim l_1l_2l_3$ since $F$ is algebraically closed, {where $e_1,e_2,e_3$ are the roots of $j_3(f_1(x,1))$, and $l_i,\ i=1,2,3$ are linear forms in $R$.} Here we repeat the discussion at the beginning of Section \ref{class-3-allchar}.

I. $l_1=l_2=l_3$.  Let $\phi(x)=l_1,\ \phi(y)=y,$ then $\phi(j_3(f_1))=x^3$, i.e. $j_3(f_1) \sim x^3.$ We have $j_3(f_1,f_2) \sim (x^3, ay^3+bxy^2+cx^2y).$

\hspace*{1em}I.1 If $a=b=c=0$, $j_3(f) \sim (x^3,0).$

\hspace*{1em}I.2 If $a=b=0,\ c \neq 0,\ j_3(f) \sim (x^3, x^2y).$

\hspace*{1em}I.3 If $a=c=0,\ b \neq 0,\ j_3(f) \sim (x^3, xy^2).$

\hspace*{1em}I.4 If $a=0,\ b,c \neq 0,\ j_3(f) \sim (x^3, bxy^2+cx^2y) \sim (x^3,xy^2+x^2y)$.

\hspace*{1em}I.5 If {$a \neq 0,\ b^2 \neq ac$}, we have $j_3(f) \sim (x^3, y^3+b'xy^2+c'x^2y) \sim (x^3,(y+b'x)^3+(c'-{b'}^2)x^2y)$, where $b'=\frac{b}{a},\ c'=\frac{c}{a}$. Using the automorphism $\phi(x)=x,\ \phi(y)=y-b'x$ we have $j_3(f) \sim (x^3,y^3+(c'-{b'}^2)x^2(y-b'x)) \sim (x^3,y^3+x^2y)$ since $b^2 \neq ac$.

\hspace*{1em}I.6 If $a \neq 0,\ b^2=3ac$, as above, we have $j_3(f) \sim (x^3,y^3+(c'-{b'}^2)x^2(y-b'x)) \sim (x^3,y^3).$

II. $l_1=l_2\neq l_3$.  Let $\phi(x)=l_1,\ \phi(y)=l_2,$ then $j_3(f_1) \sim x^2y.$ We have $j_3(f) \sim (x^2y, ax^3+bxy^2+cy^3)$

\hspace*{1em}II.1 If $a=b=c=0,\ j_3(f) \sim (x^2y,0).$

\hspace*{1em}II.2 If $a=b=0,\ c \neq 0,\ j_3(f) \sim (x^2y,y^3).$ Let $\phi(x)=y,\ \phi(y)=x,$ we get $j_3(f) \sim (x^3,xy^2).$

\hspace*{1em}II.3 If $b=c=0,\ a \neq 0,\ j_3(f) \sim (x^2y,x^3) \sim (x^3,x^2y).$

\hspace*{1em}II.4 If $a=c=0,\ b \neq 0,\ j_3(f) \sim (x^2y,xy^2).$

\hspace*{1em}II.5 If $c=0,\ a,b \neq 0,\ j_3(f) \sim (x^2y,ax^3+bxy^2)$. Let $\phi(x)=x,\ \phi(y)=\sqrt{\frac{a}{b}}y,$ we have $j_3(f) \sim (x^2y,x^3+xy^2).$

\hspace*{1em}II.6 If $c \neq 0,\ b \neq 0,\ j_3(f) \sim (x^2y,y^3+b'xy^2+a'x^3) \sim (x^2y,(y+b'x)^3+(a'-b'^3)x^3)$, where $a'=\frac{a}{c},\ b'=\frac{b}{c}$. Let $\phi(x)=x,\ \phi(y)=y-b'x,$ we have $j_3(f) \sim (x^2(y-b'x),y^3+(a'-b'^3)x^3) \sim (x^2y-b'x^3,y^3+\frac{a'-b'^3}{b'}x^2y)$. Let $\phi(x)=x,\ \phi(y)=-b'y,$ we have $j_3(f) \sim (x^3+x^2y,\ y^3+\lambda x^2y),$ where $\lambda=\frac{a'-b'^3}{b'^3}=\frac{ac^2-b^3}{b^3} \in F.$ We still assume $\lambda \neq 0$ here.

\hspace*{1em}II.7 If $c \neq 0,\ b=0,$ as above, $j_3(f) \sim (x^2(y-b'x),y^3+(a'-b'^3)x^3) \sim (x^2y, y^3+a'x^3) \sim (x^2y, x^3+y^3).$

III. $l_1\neq l_2 \neq l_3$. By multiplying a unit, one can assume $l_3=l_1+l_2.$ Suppose $l_1=x,\ l_2=y$, then $j_3(f_1) \sim x^2y+xy^2.$ Hence $j_3(f) \sim (x^2y+xy^2, ax^3+bx^2y+cy^3)$.

\hspace*{1em}III.1 If $a=b=c=0,\ j_3(f) \sim (x^2y+xy^2,0).$

\hspace*{1em}III.2 If $a=c=0,\ b \neq 0,\ j_3(f) \sim (x^2y+xy^2, x^2y) \sim (x^2y,xy^2)$.

\hspace*{1em}III.3 If exactly one of $a,c$ is equal to 0, assume $a \neq 0$. Then $j_3(f) \sim (x^2y+xy^2,x^3+b'x^2y)\sim (x^2y+xy^2,x^2(x+b'y))$, where $b'=\frac{b}{a}$. Using automorphism $\phi(x)=x,\ \phi(y)=\frac{y-x}{b'}$, it deduces to I or II depending on whether $b'=0$.

\hspace*{1em}III.4 If $a,c \neq 0$ and $b=0$, then $j_3(f) \sim (x^2y+xy^2,x^3+c'y^3) \sim (x^2y+xy^2,x^3+c'^{\frac{1}{2}}x^2y+c'^{\frac{1}{2}}xy^2+c'y^3) \sim (x^2y+xy^2,(x+cy)(x+c'^{\frac{1}{2}}y)^2)$, where $c'=\frac{c}{a}$. It's easily to {see that} it deduces to I or II depending on whether $c'=1$.

\hspace*{1em}III.5 If $a,b,c \neq 0,\ j_3(f) \sim (x^2y+xy^2,x^3+b'x^2y+c'xy^2),$ where $b'=\frac{b}{a},c'=\frac{c}{a}.$ Let $\alpha$ be the root of $\alpha^2+b'\alpha=c'$, then $f \sim (x^2y+xy^2,x^3+(\alpha+b')x^2y+\alpha xy^2+c'y^3 \sim (x^2y+xy^2,x^2(x+\lambda y)+\alpha y^2(x+\lambda y)) \sim (x^2y+xy^2,(x+\lambda y)(x+\alpha^{\frac{1}{2}}y)^2)$, where $\lambda=\alpha+b=\frac{c}{\alpha}$. Then it deduces to the case I or II.

Hence we get the result:
\begin{proposition}
    Let $f \in F[[x,y]]^2$ be a unimodal complete intersection singularity with $ord(f)=3$ in a field $F$ with characteristic 2, then $j_3(f)$ is contact equivalent to one of the following:{
    \begin{equation}
    \begin{split}
        &(x^3,0),\ (x^3,x^2y),\ (x^3,xy^2),\ (x^3,x^2y+xy^2),\ (x^3,y^3+x^2y),\ (x^3,y^3),\\
        & (x^2y,0),\ (x^2y,xy^2),\ (x^2y,x^3+xy^2),\ (x^3+x^2y,y^3+ \lambda x^2y)\ (\lambda \neq 0),\\
        &(x^2y,x^3+y^3),\ (x^2y+xy^2,0).
    \end{split}
\end{equation}}
\end{proposition}

\subsection{The classification of unimodal}
\begin{proposition}\label{ord3-class-char2}
    A unimodal ICIS of order 3 in any field with characteristic $2$ must have the form in Table \ref{table-ord3-char2}.

\renewcommand\arraystretch{1.5}
\begin{longtable}{|c|c|c|}
\caption{}\label{table-ord3-char2}
\\
\hline
Symbol&Form& condition  \\
\hline
$H_{\lambda}$&$(x^3+x^2y,y^3+\lambda x^2y)$&$\lambda \neq 0$ \\
\hline
$I$&$(x^3,y^3+x^2y)$& \\
\hline
$J$&$(x^3,y^3)$ & \\
\hline
$K_r$&$(x^3+xy^2,x^2y+ y^r)$&$r \geq 4$ \\
\hline
$L_{r,s}$&$(x^2y+y^r,xy^2+x^s)$&$r,s \geq 4$\\
\hline
$M^2$&$(x^3+y^3,x^2y)$ &\\
\hline
$M_r$&$(x^3+y^r,xy^2)$ & $r \geq 4$\\
\hline
$M_r^2$&$(x^3+y^r,x^2y+xy^2)$ &$ r \geq 4,\ r\ is\ even$\\
\hline 
$\widetilde{M}_r^2$&$(x^3+y^r+ey^l,x^2y+xy^2),\ where$ &{$ r \geq 4,\ r\ is\ odd,$}\\
&$e=e_0+e_1y^2+e_2y^4+\dots$ & $l>r,\ l\ is\ even$\\
\hline
$N_{\lambda}$&$(x^3+\lambda xy^3,x^2y+y^4)$&{$ \lambda \neq 1$}\\
\hline
$N^2_{\lambda}$&$(x^3,x^2y+y^4+\lambda xy^3)$&$ \lambda \in F$\\
\hline
$P_{r,\infty}$ & $(x^3+xy^3+xy^r,x^2y+y^4)$ & $r \geq 4$\\
\hline
$P_{\infty,s}$ & $(x^3+xy^3+y^s,x^2y+y^4)$ & $s \geq 5$ \\
\hline
$P_{r,s,\lambda}$ & $(x^3+xy^3+xy^r+\lambda y^s,x^2y+y^4)$ & $r \geq 4,s \geq 5,$ \\
&& $\lambda \in F$\\
\hline
$R_t $&$(x^3+xy^3,x^2y+y^t)$&$ t \geq 5$\\
\hline
$X^2_{\mu}$&$(x^3+y^4,x^2y+\mu xy^3)$&$\mu \in \{0,1\}$\\
\hline
$\widetilde{X}^2_{\lambda}$&$(x^3+y^4+xy^3,x^2y+\lambda xy^3)$&$\lambda \in F$\\
\hline
$Y_{\lambda}$&$(x^3+y^5,x^2y+\lambda y^5)$&$\lambda \in \{0,1\}$\\
\hline

\end{longtable}  

\end{proposition}

The following propositions in this section finish the proof.

\begin{proposition}
    If $j_3(f)$ is contact equivalent to one of the following form {
    \begin{equation}\label{*}
        (x^3,y^3+x^2y),\ (x^3,y^3),\ (x^3+x^2y,y^3+ \lambda x^2y)\ (\lambda \neq 1),\ (x^2y,x^3+y^3),
    \end{equation}}
    then $j_3(f)$ is 3-determined. In particular, $f$ is contact equivalent to the  $j_3(f)$ in {(\ref{*}).}
\end{proposition}
\begin{proof}
    For $g$ be the first three germs, we can show $\mathfrak{m}^4 \subset \mathfrak{m}\cdot \widetilde{T}_g(\mathcal{K}g)$ as before, hence $g$ is 3-determined.

    For $g=j_3(f)=(x^2y,x^3+y^3),\ \widetilde{T}_g(\mathcal{K}g)$ is spanned by {\{$(x^2y,0),(0,x^2y),(x^3+y^3,0),(0,x^3+y^3),(0,x^2),(x^2,y^2)$\}}, and a complete transversal is spanned by $\{(x^r,0)\mid r\geq 4\}$. Hence $f \sim (x^2y+e(x)x^r,x^3+y^3)$ by Theorem \ref{complete-cor} where $e(x) \in F[[x]]$ is a unit and $r\geq 4$. Applying the automorphism $\phi(x)=x,\phi(y)=y-e(x)x^{r-2},$ $f \sim (x^2y,x^3+(y+e(x)x^{r-2})^3) \sim (x^2y,x^3+y^3+e(x)^3x^{3r-6}) \sim (x^2y,u(x)x^3+y^3),$ where $u(x)=1+e(x)^3x^{3r-9}$ is a unit in $F[[x]]$. Applying $\phi(x)=x,\phi(y)=u(x)^{\frac{1}{3}}y$, then $f \sim (x^2y,x^3+y^3)$.
\end{proof}

\begin{proposition}
    If $j_3(f) \sim (x^3,xy^2)$, then {$f \sim (x^3+y^r,xy^2) \sim M_r,\ r \geq 4$ in Table \ref{table-ord3-char2}.}
\end{proposition}

\begin{proof}
    The complete transversal is still given by $\{(y^r,0) \mid r \geq 4\}$. Then the process is same as Proposition \ref{x3-xy2}.
\end{proof}

\begin{proposition}\label{x3-x2y+xy2}
    If $j_3(f) \sim (x^3,x^2y+xy^2)$, then {$f \sim M_r^2$ or $\widetilde{M}_r^2$ in Table \ref{table-ord3-char2}}. 
\end{proposition}
\begin{proof}
    A complete transversal is $T=span \langle (y^r,0), \ r \geq 4 \rangle$, then $f \sim (x^3+u(y)y^r,x^2y+xy^2)$. 
    
    If $r$ is even, using $\alpha,\beta$-trick we have $f \sim (x^3+y^r,x^2y+xy^2)$ (since $2 \nmid r-3$). 
    
    If $r$ is odd, write
    \begin{equation}
        \begin{aligned}
            f &\sim (x^3+e_0y^{2k+1}+e_1y^{2k+2}+\dots,x^2y+xy^2)\\
            &=(x^3+(e_0+e_2y^2+e_4y^4+\dots)y^{2k+1}+(e_1+e_3y^2+\dots)y^{2k+2},x^2y+xy^2).
        \end{aligned}
    \end{equation}
    There exists $e(x)^2=e_0+e_2y^2+\dots$, which allows us to use $\alpha,\beta$-trick again. Reset the symbols, {we} get a family $f \sim (x^3+y^{r}+(e_0+e_1y^2+\dots)y^l,x^2y+xy^2)$.

\end{proof}

\begin{proposition}
    If $j_3(f)\sim (x^2y,xy^2)$, then {$f \sim (x^2y+y^r,xy^2+y^s)\sim L_{r,s},\ r,s \geq 4$ in Table \ref{table-ord3-char2}.}
\end{proposition}
\begin{proof}
    The complete transversal {is the same as} Proposition \ref{rs-trick}, and the later process also follows from Proposition \ref{rs-trick}.
\end{proof}

\begin{proposition}\label{even}
    If $j_3(f) \sim (x^2y,x^3+xy^2)$, then {$f \sim (x^2y+y^r,x^3+xy^2) \sim K_r$ in Table \ref{table-ord3-char2}}, $r \geq 4$.
\end{proposition}
\begin{proof}
    Same as Proposition \ref{x2y-x3+xy2}.
\end{proof}

If {$j_3(f) \sim (x^3,x^2y)$}, a complete transversal is given by $$\{(xy^r,0),(y^s,0),(0,xy^u),(0,y^v) \mid r,u \geq3,s,v \geq 4\}.$$ {\color{black}Then} 
\begin{equation}\label{char2-x3-x2y}
    f \sim (x^3+a(y)xy^r+b(y)y^s,x^2y+c(y)xy^u+d(y)y^v).
\end{equation}
with {$r,u \geq 3,\ s,v \geq 4$}

\begin{proposition}\label{bimod-ord3-char2}
    If $f$ is of the form {$(\ref{char2-x3-x2y})$} and $r\geq 4,\ s \geq 6,\ v \geq 5$, then $f$ is at least bimodular. 
\end{proposition}
\begin{proof}
    For {$r \geq 4,\ s\geq 6,v \geq 5$, we write $g=j_4(f)=(x^3,x^2y)$, and any 5-jet in a open dense subset of $J'_5(g)$ is of the form $g_{ac}=(x^3+axy^4,x^2y+cxy^3+dy^5)$ with $a,c,d\in F$, where $J'_5(g)$ is formed by jets in $J_5(g)$ with $r \geq 4,s \geq 6,v \geq 5$.} Computation or Singular program can show $(y^5,0),(0,y^5) \notin P_{4,5}$ in each case, hence $cod(g_{ac})=2$ for all $a,c \in F$. By Proposition \ref{cod-modality} and \ref{modality}, $\mathcal{K}-mod(f) \geq \mathcal{K}_5-mod(f) \geq inf\{cod(g_{ac})\} \geq 2.$
\end{proof}

\begin{proposition}
    If $f$ is of the form {$(\ref{char2-x3-x2y})$} with $s=4$, then $f$ is contact equivalent to the form $\sim (x^3+y^4,x^2y+\mu xy^3)$ with $\mu=0,1$ or $(x^3+y^4+xy^3,x^2y+\lambda xy^3)$ with $\lambda \in F$. That is, {$f \sim X^2_{\mu}$ or $\widetilde{X}^2_{\lambda}$ in Table \ref{table-ord3-char2}.}
\end{proposition}
\begin{proof}
    In this case $f \sim (x^3+y^4+a(y)xy^r,x^2y+c(y)xy^u+d(y)y^v)$ after a scaler transform. Let $g=(x^3+y^4,x^2y)$, then $g$ is weighted homogeneous of degree 0 with respect to $(a,d)$, where $a=(4,3),\ d=(12,11)$. Let $T=span\langle (xy^3,0),(0,xy^3) \rangle$, then we can check that {$F_{a,d}^1R^2 \subset T+\widetilde{T}_g(F^1\mathcal{K}g)$}, which means $f \sim (x^3+y^4+axy^3,x^2y+bxy^3)$ for $a,b \in F$. Take a scaler transform, we have $f \sim (x^3+y^4,x^2y)$, $(x^3+y^4,x^2y+xy^3)$ or $(x^3+y^4+xy^3,x^2y+\lambda xy^3)$ with $\lambda \in F$.
\end{proof}

\begin{proposition}
    If $f$ is of the form {$(\ref{char2-x3-x2y})$} with $v=4,s > 4$, then {$f \sim N_{\lambda},N^2_{\lambda}$ or $ P_{r,\infty}, P_{\infty,s}, P_{r,s,\lambda}$ in Table \ref{table-ord3-char2}.}
\end{proposition}
\begin{proof}
    If $r \geq 4,$ then let $g=(x^3,x^2y+y^4)$ be the weighted 0-jet of $f$ with respect to $(a,d)$, where $a=(3,2),d=(9,8)$. Let $T=span\langle(0,xy^3) \rangle$, then {$F_{a,d}^1R^2 \subset T+\widetilde{T}_g(F^1\mathcal{K}g)$}, which means $f \sim (x^3,x^2y+y^4+\lambda xy^3)$ for $\lambda \in F$.

    If $r=3$, then let $g=(x^3+a(0)xy^3,x^2y+y^4),\ a(0)\in F^{\times}$ be the weighted 0-jet of $f$ with respect to $(a,d)$, where $a=(3,2),d=(9,8)$. If $a(0) \neq 1$, then {$F_{a,d}^1R^2 \subset \widetilde{T}_g(F^1\mathcal{K}g)$} automatically holds, hence $f \sim (x^3+\lambda xy^3,x^2y+y^4),\ \lambda \in F^{\times}$.
    
    In the case $a(0)=1$, $T$ is spanned by $\{(xy^r,0),(y^s,0)\mid r \geq 4, s \geq 5\}$, and $f \sim (x^3+xy^3+u(y)xy^r+v(y)y^s,x^2y+y^4)$. Similar as Proposition \ref{s6t4}, for $v(y)=0$ (resp. $u(y)=0$), we have $f \sim P_{r,\infty}$ (resp. $P_{\infty,s}$). Otherwise, since $p =2$, $p \nmid 2r-2s-3$ holds for any $r,s$. Then we have $f \sim (x^3+xy^3+e(y)xy^r+e(y)y^s,x^2y+y^4) \sim (x^3+xy^3+e_0xy^r+e_0y^s,x^2y+y^4) \sim (x^3+xy^3+xy^r+\lambda y^s,x^2y+y^4) \sim P_{r,s,\lambda}$.
\end{proof}

\begin{proposition}
    If $f$ is of the form {$(\ref{char2-x3-x2y})$} with $s=5$, then {$f \sim Y_{\lambda}$ or $R_t$ in Table \ref{table-ord3-char2}.}
\end{proposition}
\begin{proof}
    If $r \geq 4,$ then let $g=(x^3+by^5,x^2y)$ with $b \in F^{\times}$ be the weighted 0-jet of $f$ with respect to $(a,d)$, where $a=(5,3),d=(15,13)$. Same to Proposition \ref{s5r4}, $f \sim (x^3+y^5,x^2y+\lambda y^5),\ \lambda \in \{0,1\}$.

    If $r=3$, let $g=(x^3+axy^3,x^2y)$ with $a \in F^{\times}$ be the weighted 0-jet of $f$ with respect to $(a,d)$, where $a=(3,2),d=(9,8)$. Same to Proposition \ref{s-geq-6}, we have $f \sim (x^3+xy^3,x^2y+y^t),\ t \geq 5$.
\end{proof}

\begin{proposition}
    If $f$ is of the form {$(\ref{char2-x3-x2y})$} with $s \geq 6,v \geq 5$, then {$f \sim R_t$ in Table \ref{table-ord3-char2}.}
\end{proposition}
\begin{proof}
    In this case we have $r=3$ by Proposition \ref{bimod-ord3-char2}. Under the assumption $s,v>4$, we can choose $g=(x^3+axy^3,x^2y)$ be the weighted 0-jet of $f$ with respect to $(a,d)$, where $a=(3,2),d=(9,8)$. Let $T=span\langle(0,y^l)|\ l \geq 5 \rangle$, then {$F_{a,d}^1R^2 \subset T+\widetilde{T}_g(F^1\mathcal{K}g)$}. Therefore $f \sim (x^3+xy^3,x^2y+e(y)y^t)$. Using $r,s$-trick same as Proposition \ref{s-geq-6}, we have $f \sim (x^3+xy^3,x^2y+y^t),\ t \geq 5$.
\end{proof}

\section{The classification of order 3 when $char F=3$}
Next we repeat the discussion in the case $charF=3$.

\subsection{The classification of {3-jets}}
Same as in Section \ref{class-3-allchar}, for $f=(f_1,f_2)$ with $\mathrm{ord}(f_1)=3$, we have {$j_3(f_1) \sim x^3,x^2y$ or $x^3+xy^2$}.

I. $j_3(f_1) \sim x^3.$ We have $j_3(f_1,f_2) \sim (x^3, ay^3+bxy^2+cx^2y).$

\hspace*{1em}I.1 If $c = 0,\ j_3(f) \sim (x^3, x^2y)$ or $(x^3,xy^2)$ or $(x^3,y^3+xy^2)$ or $(x^3,0)$ depending on whether $a,b=0$.

\hspace*{1em}I.2 If $c \neq 0,\ b =0,\ j_3(f) \sim (x^3,y^3+x^2y).$

\hspace*{1em}I.3 If $b,c \neq 0,$ then $j_3(f) \sim (x^3,y(x+\frac{b}{2}y)^2+(a-\frac{b^2}{4})y^3) \sim (x^3,\frac{b^2}{4}y(y+\frac{2}{b}x)^2+(a-\frac{b^2}{4})(y+\frac{2}{b}x)^3).$ Using $\phi(x)=x,\phi(y)=y+\frac{2}{b}x,\ j_3(f) \sim (x^3,\frac{b^2}{4}(y-\frac{2}{b}x)y^2+(a-\frac{b^2}{4})y^3)\sim (x^3,ay^3-\frac{b}{2}xy^2)\sim (x^3,xy^2)$ or $(x^3,y^3+xy^2)$ depending on whether $a=0$.

II. {$j_3(f_1) \sim (x^2y,0).$} We have {$j_3(f) \sim (x^2y, ax^3+bxy^2+cy^3)$.}

\hspace*{1em}II.1 If $c=0,\ j_3(f) \sim (x^2y,ax^3+bxy^2) \sim (x^2y,x^3)$ or $(x^2y,xy^2)$ or $(x^2y,x^3+xy^2)$ or $(x^2y,0)$ depending on whether $a,b$ equal to 0. 

\hspace*{1em}II.2 If $c \neq 0,\ b=0,\ j_3(f) \sim (x^2y,ax^3+cy^3) \sim (x^2,(a^{\frac{1}{3}}x+c^{\frac{1}{3}}y)^3),$ which is back to case I.

\hspace*{1em}II.3 If $b,c \neq 0,\ a=0,\ j_3(f) \sim (x^2y,y^3+xy^2)$.

\hspace*{1em}II.4 If $a,b,c \neq 0,$ change the notations, $ j_3(f) \sim (x^2y,y^3+ax^3+bxy^2).$ {Applying the automorphism $\phi(x)=x,\phi(y)=y+\alpha x$, where $\alpha$ is a nonzero root of $\alpha^3-b\alpha^2+a=0$ (in this case $\alpha \neq b$ since $a \neq 0$), we have $j_3(f) \sim (x^2y+\alpha x^3,y^3+(\alpha^3+b\alpha^2+a)x^3+bxy^2+2b\alpha x^2y) \sim (x^2y+\alpha x^3,y^3+bxy^2+2b\alpha(x^2y+\alpha x^3)) \sim (x^2y+\alpha x^3,y^3+bxy^2) \sim (x^3+x^2y,y^3+\lambda xy^2),\lambda \neq 1$ since $\alpha \neq b$.}

III. $j_3(f_1) \sim x^3+xy^2.$ We have $j_3(f) \sim (x^3+xy^2, axy^2+bx^2y+cy^3)$. All the {discussions} is {same as that} in Section \ref{class-3-allchar}, except {for the characteristic becomes to $3$}, with nothing changes.

In conclusion, we get:

\begin{proposition}
    Let $f \in F[[x,y]]^2$ be a unimodal complete intersection singularity with $ord(f)=3$ in a field $F$ with characteristic 3, then $j_3(f)$ is contact equivalent to one of the following:
    \begin{equation}
    \begin{split}
        &(x^3,0),\ (x^3,x^2y),\ (x^3,xy^2),\ (x^3,y^3+xy^2),\ (x^3,y^3+x^2y),\ (x^3,y^3),\\
        & (x^2y,0),\ (x^2y,xy^2),\ (x^2y,x^3+xy^2),\ (x^2y,y^3+xy^2),\ (x^3+x^2y,y^3+ \lambda xy^2)(\lambda \neq 1),\\
        & (x^3+xy^2,0).
    \end{split}
\end{equation}
\end{proposition}

\subsection{The classification of unimodal}
\begin{proposition}\label{ord3-class-char3}
    A unimodal ICIS of order 3 in any field with characteristic $3$ must have the form in Table \ref{table-ord3-char3}.
\renewcommand\arraystretch{1.5}
\begin{longtable}{|c|c|c|}
\caption{}\label{table-ord3-char3}
\\
\hline
Symbol&Form& condition  \\
\hline
$H$&$(x^3+x^2y,y^3+\lambda x^2y)$& {$\lambda \neq 0$} \\
\hline
$I^3_{\lambda}$&$(x^3+\lambda y^4,y^3+x^2y)$&$ \lambda \in \{0,1\}$\\
\hline
$\widetilde{I}^3_{\lambda}$&$(x^3+\lambda y^4,y^3+xy^2)$&$ \lambda \in \{0,1\}$\\
\hline
$J^3_{\lambda,\mu}$&$(x^3+\lambda x^2y^2,y^3+\mu x^2y^2)$ & $\lambda,\mu \in \{0,1\}$\\
\hline
$K_r$&$(x^3+xy^2,x^2y+ y^r)$&$r \geq 4$ \\
\hline
$K^3_r$&$(y^3+xy^2+x^r,x^2y)$ & $r \geq 4$\\
\hline
$\widetilde{K}^3_{r,\lambda}$&$(y^3+xy^2+\lambda xy^3+x^r,x^2y)$ & $r \geq 4,\lambda \in F,3 \mid r $\\
\hline
$L_{r,s}$&$(x^2y+y^r,xy^2+x^s)$&$r,s \geq 4$\\
\hline
$M^3_r$&$(x^3+y^r,xy^2)$ & $r \geq 4$\\
\hline
$\widetilde{M}^3_{r,\lambda}$&$(x^3+\lambda x^3y+y^r,xy^2)$ & $r \geq 4,\lambda \in F,3 \mid r$\\
\hline
$N^3_{\lambda}$&$(x^3+\lambda xy^3,x^2y+y^4)$&$ \lambda \neq 1$\\
\hline
$\widetilde{N}^3_{s}$&$(x^3+y^s,x^2y+y^4)$&$ s \in \{5,6\}$\\
\hline
$P_{r,\infty}$ & $(x^3+xy^3+xy^r,x^2y+y^4)$ & $r \geq 4$\\
\hline
$P_{\infty,s}$ & $(x^3+xy^3+y^s,x^2y+y^4)$ & $s \geq 5$ \\
\hline
$P_{r,s,\lambda}$ & $(x^3+xy^3+xy^r+\lambda y^s,x^2y+y^4)$ & $r \geq 4,s \geq 5,$ \\
&& $\lambda \in F$\\
\hline
$\widetilde{P}_{r,s} $&$(x^3+xy^3+uxy^r+vy^s,x^2y+y^4),\ where$& $r \geq 4, \ s \geq 5,$\\
&{$u=u_0+u_1y+\dots,\ v=v_0+v_1y+\dots,$}& {$3 \mid 2r-2s$}\\
\hline
$R_t$ &$(x^3+xy^3,x^2y+y^t)$& $t \geq 5$\\
\hline
$X_{\lambda}$&$(x^3+y^4,x^2y+\lambda y^4)$&$\lambda \in \{0,1\}$\\
\hline
$Y_{\lambda}$&$(x^3+y^5,x^2y+\lambda y^5)$&$\lambda \in \{0,1\}$\\
\hline

\end{longtable}  
\end{proposition}

{\color{black}Since most of the discussion is the same as before}, we make a table to present the results. See Table \ref{table char3-process}.

\renewcommand\arraystretch{1.5}
\begin{longtable}{|c|c|c|}
\caption{}\label{table char3-process}
\\
\hline
$j_3(f)$&complete\ transversal &$f$\\
\hline
$(x^3,xy^2)$&$(y^r,0), r \geq 4$&$(x^3+e(y)y^r,xy^2),\ r\geq 4$\\
\hline
$(x^3,y^3+xy^2)$&$(y^4,0)$&$(x^3+\lambda y^4,y^3+xy^2),\lambda \in \{0,1\}$\\
\hline
$(x^3,y^3+x^2y)$&$(y^4,0)$&$(x^3+\lambda y^4,y^3+x^2y),\lambda \in \{0,1\}$\\
\hline
$(x^3,y^3)$&$(x^2y^2,0),(0,x^2y^2)$&$(x^3+\lambda x^2y^2,y^3+\mu x^2y^2),\lambda,\mu \in \{0,1\}$\\
\hline
$(x^2y,xy^2)$&$(y^r,0),(0,x^s),r,s \geq 4$&{$(x^2y+y^r,xy^2+x^s),r,s \geq 4$}\\
\hline
$(x^2y,x^3+xy^2)$&$(y^r,0),r \geq 4 $&$(x^2y+y^r,x^3+xy^2),\ r \geq 4$\\
\hline
$(x^2y,y^3+xy^2)$&$(0,x^r),r \geq 4 $&$(x^2y,y^3+xy^2+e(x)x^r),\ r \geq 4$\\
\hline
$(x^3+x^2y,$&$ 3-determined $&$(x^3+x^2y,y^3+\lambda xy^2)$\\
$y^3+\lambda xy^2)$&&\\
\hline
$(x^3,x^2y)$&$(y^r,0),(xy^s,0),(0,y^t)$&$(x^3+a(y)y^r+b(y)y^s,xy^2+c(y)y^t)$\\
\hline

\end{longtable}

When $j_3(f) \sim (x^3,xy^2),(x^2y,y^3+xy^2),(x^3,x^2y)$, we need further {discussion}.

\begin{proposition}
    If $j_3(f) \sim (x^3,xy^2)$, then $f \sim M^3_r$ or $ \widetilde{M}^3_{r,\lambda}$.
\end{proposition}
\begin{proof}
    We have $f \sim (x^3+e(y)y^r,xy^2)$ with $r \geq 4,\ e(y) \in F[[y]]$ is a unit.

    If $3 \nmid r$, using $\alpha,\beta$-trick we have $f \sim (x^3+y^r,xy^2)$.

    If $3 \mid r$, reset notations, write $f \sim (e(y)x^3+y^r,xy^2) \sim (x^3+e_1x^3y+y^r,xy^2)$, where {$e_1 \in F$}. Then $f \sim \widetilde{M}^3_{r,e_1}$.

\end{proof}

\begin{proposition}
    If $j_3(f) \sim (x^2y,y^3+xy^2)$, then $f \sim K^3_r$ or $\widetilde{K}^3_{r,\lambda}$.
\end{proposition}
\begin{proof}
    We have $f \sim (x^2y,y^3+xy^2+e(x)x^r).$ 
    
    If $3 \nmid r$, using $\alpha,\beta$-trick we have $f \sim (x^2y,y^3+xy^2+x^r)\sim K^3_r.$
    
    If $3 \mid r$, similar as above, $f \sim (y^3+xy^2+\lambda xy^3+x^r,x^2y) \sim \widetilde{K}^3_{r,\lambda}$ with $\lambda \in F$.
\end{proof}

When $j_3(f) \sim (x^3,x^2y)$, the following result is the same as the case $charF>3$.
\begin{proposition}
    If $j_3(f) \sim (x^3,x^2y)$, then $f \sim (x^3+a(y)y^r+b(y)y^s,xy^2+c(y)y^t)$ with $r,t\geq 4,\ s \geq 3$, $a(y),b(y),c(y)$ are units. And when $r \geq 4,\ s \geq 6,\ t \geq 5$, $f$ is not unimodal.
\end{proposition}

Here we omit the discussion same as Proposition \ref{s4} to {Proposition} \ref{s-geq-6} and give the result directly in Table \ref{eqtable}.

\begin{longtable}{|c|c|c|c|c|}
\caption{}\label{eqtable}
\\
\hline
&weighted\ jet&weight &complete& form \\
&&&transversal&\\
\hline
$s=4$&$(x^3+y^4,x^2y)$&$(4,3; $&$ (0,y^4) $&$ (x^3+y^4,x^2y+\mu y^4),$\\
&&$12,11)$&&$\mu \in \{0,1\}$\\
\hline
$s=5,$&$(x^3+xy^3,x^2y+\lambda y^4)$&$(3,2;9,8)$&$ 0 $&{$(x^3+xy^3,x^2y+\lambda y^4),$}\\
$r=3$&$\lambda \notin \{0,1\}$&&&{$\lambda \notin \{0,1\}$}\\
\hline
$s=5,$&$(x^3+xy^3,x^2y+\lambda y^4) $&$(3,2;9,8)$&$ (0,y^t),t \geq 5 $&{$(x^3+xy^3,x^2y+ y^t),$}\\
$r=3$&$\lambda =0$&&&{$t \geq 5 $}\\
\hline
$s=5,$&$(x^3+xy^3,x^2y+\lambda y^4)$&$(3,2;9,8)$&$ {(xy^r,0),(y^s,0)} $&${(x^3+xy^3+u(y)xy^r}$ \\
$r=3$&$\lambda=1$&&& ${+v(y)y^s,x^2y+y^4)}$\\
&&&&{where $u(y),v(y)$ are }\\
&&&&{units in $F[[x,y]]$,}\\
&&&&{$r \geq4,s\geq5$}\\
\hline
$s=5,$&$(x^3,x^2y+y^4)$&$(3,2;9,8)$&$ (y^5,0),(y^6,0) $&$(x^3+\lambda y^s,x^2y+y^4),$\\
$r \geq 4,$&&&&$s=5,6,\lambda=0,1$\\
$t=4$&&&&\\
\hline
$s=5,$&$(x^3+y^5,x^2y)$&$(5,3; $&$ (0,y^5) $&$ (x^3+y^5,x^2y+\mu y^5),$\\
$r \geq 4,$&&$15,13)$&&$\mu \in \{0,1\}$\\
$t \geq 5$&&&&\\
\hline
$s \geq 6, $&$(x^3+\lambda xy^3,x^2y+y^4)$&$(3,2;9,8)$&$0 $&$ (x^3+\lambda xy^3,x^2y+y^4),$\\
$t=4$&$\lambda \notin \{0,1\}$&&&$\lambda \notin \{0,1\}$\\
\hline
$s \geq 6, $&$(x^3+\lambda xy^3,x^2y+y^4)$&$(3,2;9,8)$&$ (y^5,0),(y^6,0)$&$ (x^3+\lambda y^s,x^2y+y^4),$\\
$t=4$&$\lambda=0$&&&$s=5,6,\lambda=0,1$\\
\hline
$s \geq 6, $&$(x^3+\lambda xy^3,x^2y+y^4)$&$(3,2;9,8)$&$(xy^r,0),(y^s,0)$&${(x^3+xy^3+u(y)xy^r}$\\
$t=4$&&&&${+v(y)y^s,x^2y+y^4)}$\\
&&&&{where $u(y),v(y)$ are }\\
&&&&{units in $F[[x,y]]$,}\\
&&&&{$r \geq4,s\geq5$}\\
\hline
$r=3,$&$(x^3+xy^3,x^2y)$&$(3,2;9,8)$&$ (0,y^t),t \geq 5 $&$(x^3+xy^3,x^2y+ y^t),$\\
$s\geq 6,$&&&&{$t \geq 5$}\\
$t\geq 5$&&&&\\
\hline

\end{longtable}

These finish the proof of Proposition \ref{ord3-class-char3}.

\section{Check the modality}
Let $T^{1,sec}(f)=R^m \bigg/ \bigg(\langle f_1,\dots,f_m \rangle \cdot R^m+\mathfrak{m} \langle \frac{\partial f}{\partial x_1},\dots,\frac{\partial f}{\partial x_n} \rangle \bigg)$. Choose $g_1,\dots,g_d$ to be a $F$-basis of $T^{1,sec}(f)$. T.H. Pham, G. Pfister, and G.M. Greuel {have} shown in \cite{icissimple} that 
{$$F_{\mathbf{t}}(x)=F(x,t_1,\dots,t_d)=f+t_1g_1+\dots+t_dg_d $$ 
represents a formally semiuniversal deformation of $f$, where $\mathbf{t}=(t_1,\dots,t_d)$.} If $F(x,t)$ is equivalent to a family of ICIS of at most one parameter for $t \in F^d$, then $f$ is unimodal. Here we check $f \sim l_{q,\lambda} \sim (x^2+y^4,y^{q+2}+\lambda xy^q),\ q \geq 3,\lambda^2 \notin \{0,-1\}$ in Table \ref{table-ord2} as an {example.}

First we choose generators $$(y,0),(y^2,0),(y^3,0),(0,y),(0,y^2),\dots,(0,y^{q+2}) \in T^{1,sec}(f).$$ Note that {$T^{1,sec}(f)=(\mathfrak{m}R^2)/\widetilde{T}_f(\mathcal{K}f)$.} In the proof of {\textbf{(3.b.ii.2.3.3.1)} in Proposition \ref{ord2-class-mod1}}, we have shown that $(0,xy^{q+1}) \in \widetilde{T}_f(\mathcal{K}f)$. That means $(0,y^{q+3}) \in \widetilde{T}_f(\mathcal{K}f)$. We also have $(0,y^{q+2})$ generates $(0,xy^{q})$, $(y^3,0)$ and $(0,y^{q+1})$ generate $(0,xy^{q-1})$. Then we add $(0,x),(0,xy),\dots,(0,xy^{q-2})$ as generators. These generators form a basis of $T^{1,sec}(f)$.

Let $$g_1=(y^2,0),g_2=(y^3,0),g_3=(0,y^2),\dots, g_{q+3}=(0,y^{q+2}),$$ $$g_{q+4}=(0,xy),\dots, g_{2q+1}=(0,xy^{q-2}).$$ {We consider $F(x,\mathbf{t})=f+t_1g_1+\dots+t_{2q+1}g_{2q+1}$, where $\mathbf{t}=(t_1,\dots, t_{2q+1})$. We write $F(x,\mathbf{t})=(G_1,G_2) \in R^2$. }

If $t_1 \neq 0$ or $t_3 \neq 0$ or $t_{q+4} \neq 0$, then $j_2(G_1)$ is non-degenerate, which means $G$ is simple by Proposition \ref{simple-form}. From now we assume $t_1=t_3=t_{q+4}=0$.

If $t_2 \neq 0$ or $t_{4} \neq 0$, then $j_2(G_1) \sim (x^2+y^3)$. By Proposition \ref{simple-form}, $G$ is simple. From now we assume $t_2=t_{4}=0$.

Now $G$ is of the form 
\begin{equation}
    \begin{aligned}
        G &\sim (x^2+y^4,y^{q+2}+\lambda y^{q}+t_5y^4+\dots+t_{q+3}y^{q+2}+t_{q+5}xy^2+\dots+t_{2q+1}xy^{q-2})\\
        & \sim(x^2+y^4,\sum_{i\geq u}t_{i+1}y^i+x\sum_{j\geq v}t_{j+q+3}y^{j}),
    \end{aligned}
\end{equation}{where $q \geq 3,\lambda^2 \notin \{0,1\},u\geq 4,v\geq 2$.}
Comparing with {(\ref{ord2form})} and Proposition \ref{ord2-class-mod1}, it corresponds to the case $\alpha=1,\ s=4,\ u=t\geq 4,\ v=q \geq 2$, where $\alpha,s,t,q$ is in the sense of {(\ref{ord2form}).}

By Proposition \ref{ord2-class-mod1} {\textbf{(3.a.i)}}, in most of cases, $G$ is at most unimodal. The only unsure case is that there exists $q'<t<t'$ such that 
$$G \sim (x^2+y^4,y^{q'+2}+\lambda' xy^{q'}+u_0xy^{t}+xy^{t'}+u_1xy^{t+p})$$ 
with $q' \geq 3$, $t+p\leq q-2$, ${\lambda'}^2=-1$ and $p \mid t-q',\ p \nmid t'-q'$, which is of the form $\widetilde{l}_{q',t,t'}$ {in Table \ref{table-ord2}} containing two parameters $u_0$ and $u_1$. Then we must have $$p+3 \leq p+q' \leq t  \leq q-p-2,$$ that is, $q \geq 2p+5$. If we set $q \leq 2p+4,$ then this will not happen. That is, $l_{q,\lambda}$ is unimodal for $q \leq 2p+4$.

Use the same method, we can give tables to show when the modality of the above class is $1$.

\renewcommand\arraystretch{1.5}
\begin{longtable}{|c|c|c|c|}
\caption{}\label{table-ord2-mod1}
\\
\hline
Symbol&Form& condition&when is  \\
&&&unimodal\\
\hline
$h_q$&$(x^2+y^4,xy^q)$&$q \geq 3$ & $q \leq 2p+3$ \\
\hline
$i$&$(x^2,y^5)$&& \\
\hline
$\widetilde{i}$&$(x^2,y^5+xy^3)$&& \\
\hline
$i^5$&$(x^2,y^5+xy^4)$& $p=5$ &\\
\hline
$j_t$&$(x^2+y^4,y^t)$ &$t \geq 5$ & $t \leq 2p+4$ \\
\hline
$\widetilde{j}_t$&$(x^2+y^4,y^t+xy^{t-1})$&$t \geq 5,\ p \mid t$ & $t\leq 2p+4$ \\
\hline
$k_q$&$(x^2+y^4,y^{q+3}+xy^q)$&$q \geq 3$ & $q \leq 2p+3$\\
\hline
$l_{q,\lambda}$&$(x^2+y^4,y^{q+2}+\lambda xy^q)$ & $q \geq 3, \lambda^2 \notin \{0, -1\}$ & $q \leq 2p+4$\\
\hline
$\widetilde{l}_{q,t,t'} $&{$(x^2+y^4,y^{q+2}+\lambda xy^q+uxy^t+xy^{t'}),$} & $\lambda^2=-1, q \geq 3,$ & $q \leq 2p+3 $ \\
&{where $u=u_0+u_1y^p+u_2y^{2p}+\dots$}& $ t \geq q+1, t' \geq t+1, $& \\
&& $p \mid t-q, p \nmid t'-q$&\\
\hline
\end{longtable}

\begin{remark}
    If we have {that} $$(x^2+y^4,y^{5}+\lambda xy^3+axy^{p+3}+xy^{4+p}+bxy^{2p+3}) $$ is not contact equivalent to $$(x^2+y^4,y^{5}+\lambda xy^3+cxy^{p+3}+xy^{4+p}+dxy^{2p+3}) $$ for general $a,b,c,d \in F$, then we can ensure that the singularities given in Table \ref{table-ord2-mod1} are the only unimodal ICIS. Conversely, if all singularities of the form $\widetilde{l}_{q,t,t'}$ are equivalent (or at least can be presented as {a one-parameter} family), then all the singularities given in Table \ref{table-ord2} are unimodal.
    Unfortunately, we cannot judge this equivalence yet. We post it as a conjecture.
\end{remark}

\begin{conjecture}\label{conjecture}
    Let $F$ be an algebraically closed field with characteristic $p$. Then the isolated complete intersection singularity
    $$(x^2+y^4,y^{5}+\lambda xy^3+axy^{p+3}+xy^{4+p}+bxy^{2p+3}) $$ is not contact equivalent to $$(x^2+y^4,y^{5}+\lambda xy^3+cxy^{p+3}+xy^{4+p}+dxy^{2p+3}) $$ for general $a,b,c,d \in F$.
\end{conjecture}

The modality of singularities in Table \ref{table-ord2-char2} (i.e. of order $2$ in characteristic $2$ field) do not need to check, since every germs of the form $ (x^2+h,g)$ with $g \in \mathfrak{m}^4 \backslash \mathfrak{m}^5$ is equivalent to a form in Table \ref{table-ord2-char2}, which has at most one parameter.

Using the same method to check Table \ref{table-ord3}, we find that:

\renewcommand\arraystretch{1.5}
\begin{longtable}{|c|c|c|c|}
\caption{}\label{table-ord3-mod1}
\\
\hline
Symbol&Form& condition&when is  \\
&&&unimodal\\
\hline
{$H$}&$(x^3+x^2y,y^3+\lambda x^2y)$&{$\lambda \neq 0$} &\\
\hline
$I$&$(x^3,y^3+x^2y)$&& \\
\hline
$J$&$(x^3,y^3)$ && \\
\hline
$K_r$&$(x^3+xy^2,x^2y+ y^r)$&$r \geq 4$ & \\
\hline
$L_{r,s}$&$(x^2y+y^r,xy^2+x^s)$&$r,s \geq 4$ &\\
\hline
$M_r$&$(x^3+y^r,xy^2)$ & $r \geq 3$ &\\
\hline
$N_{\lambda}$&$(x^3+\lambda xy^3,x^2y+y^4)$&$ \lambda \notin \{1,12\}$ &\\
\hline
$P_{r,\infty}$ & $(x^3+xy^3+xy^r,x^2y+y^4)$ & $r \geq 4$&\\
\hline
$P_{\infty,s}$ & $(x^3+xy^3+y^s,x^2y+y^4)$ & $s \geq 5$ &\\
\hline
$P_{r,s,\lambda}$ & $(x^3+xy^3+xy^r+\lambda y^s,x^2y+y^4)$ & $r \geq 4,s \geq 5,$& \\
&& $\lambda \in F$&\\
\hline
$\widetilde{P}_{r,s} $&$(x^3+xy^3+uxy^r+vy^s,x^2y+y^4),\ where$& $r \geq 4, \ s \geq 5$ &\\
&$u=u_0+u_1y+\dots,\ v=v_0+v_1y+\dots,$& $p \mid 2r-2s+3$ &\\
\hline
$R_t$ &$(x^3+xy^3,x^2y+y^t)$&$ t \geq 5$ &\\
\hline
$X_{\lambda}$&$(x^3+y^4,x^2y+\lambda y^4)$&$\lambda \in \{0,1\}$&\\
\hline
$Y_{\lambda}$&$(x^3+y^5,x^2y+\lambda y^5)$&$\lambda \in \{0,1\}$&\\
\hline
$Z_{\lambda}$&$(x^3+12xy^3+\lambda y^5,x^2y+y^4)$&$\lambda \in \{0,1\}$&\\
\hline
\end{longtable}

That means, every singularities in Table \ref{table-ord3} is unimodal. This is because the bimodular one must have a deformation to $\widetilde{P}_{r,s}$. But $(xy^k,0),(y^l,0)\ (k \geq 4, l \geq 5)$ do not exist in a basis of $T^{1,sec}(f)$ at the same time for every $f$. Then such a deformation {does} not exist.

For $ord(f)=3,\ char F=2$ case, see Table \ref{table-ord3-char2-mod1}:
\renewcommand\arraystretch{1.5}
\begin{longtable}{|c|c|c|c|}
\caption{}\label{table-ord3-char2-mod1}
\\
\hline
Symbol&Form& condition&when is  \\
&&&unimodal\\
\hline
{$H$}&$(x^3+x^2y,y^3+\lambda x^2y)$&$\lambda \neq 0$& \\
\hline
$I$&$(x^3,y^3+x^2y)$& &\\
\hline
$J$&$(x^3,y^3)$ && \\
\hline
$K_r$&$(x^3+xy^2,x^2y+ y^r)$&$r \geq 4$& \\
\hline
$L_{r,s}$&$(x^2y+y^r,xy^2+x^s)$&$r,s \geq 4$&\\
\hline
$M^2$&$(x^3+y^3,x^2y)$ &&\\
\hline
$M_r$&$(x^3+y^r,xy^2)$ & {$r \geq 4$}& $r \leq 8$\\
\hline
$M_r^2$&$(x^3+y^r,x^2y+xy^2)$ &$ r \geq 4,\ r\ is\ even$& $r \leq 8$\\
\hline
$\widetilde{M}_r^2$&$(x^3+y^r+ey^l,x^2y+xy^2),$ &{$ r \geq 4,\ r\ is\ odd,$}& $l \leq 7$\\
&where $e=e_0+e_1y^2+e_2y^4+\dots$ & {$l>r,\ l\ is\ even$}&\\
\hline
$N_{\lambda}$&$(x^3+\lambda xy^3,x^2y+y^4)$&$ \lambda \neq 1$&\\
\hline
$N^2_{\lambda}$&$(x^3,x^2y+y^4+\lambda xy^3)$&$ \lambda \in F$&\\
\hline
$P_{r,\infty}$ & $(x^3+xy^3+xy^r,x^2y+y^4)$ & $r \geq 4$& $r \leq 7$\\
\hline
$P_{\infty,s}$ & $(x^3+xy^3+y^s,x^2y+y^4)$ & $s \geq 5$ &$s \leq 8$ \\
\hline
$P_{r,s,\lambda}$ & $(x^3+xy^3+xy^r+\lambda y^s,x^2y+y^4)$ & $r \geq 4,s \geq 5,$& If $r \leq s+1,r \leq 7$ \\
&& $\lambda \in F$ & If $s \leq r, s \leq 8$\\
\hline
$R_t $&$(x^3+xy^3,x^2y+y^t)$&$ t \geq 5$&\\
\hline
$X^2_{\mu}$&$(x^3+y^4,x^2y+\mu xy^3)$&$\mu \in \{0,1\}$&\\
\hline
$\widetilde{X}^2_{\lambda}$&$(x^3+y^4+xy^3,x^2y+\lambda xy^3)$&$\lambda \in F$&\\
\hline
$Y_{\lambda}$&$(x^3+y^5,x^2y+\lambda y^5)$&$\lambda \in \{0,1\}$&\\
\hline
\end{longtable}  

When an ICIS in Table \ref{table-ord3-char2-mod1} can deform to $\widetilde{M}^2_r$, it may have two parameters. This leads to another conjecture:
\begin{conjecture}\label{conjecture2}
    Let $F$ be an algebraically closed field with characteristic $2$. Then the isolated complete intersection singularity
    $$(x^3+y^5+ay^6+by^8,x^2y+xy^2) $$ is not contact equivalent to $$(x^3+y^5+cy^6+dy^8,x^2y+xy^2) $$ for general $a,b,c,d \in F$.
\end{conjecture}
If the conjecture holds, then all unimodal ICIS of order $3$ in a characteristic $2$ field are presented in Table \ref{table-ord3-char2-mod1}.

For $ord(f)=3,\ char(F)=3$ case, we can check that every singularities in Table \ref{table-ord3-char3} is unimodal.

In conclusion, we get the following classification theorem:
\begin{theorem}\label{main-result}
    Let $F$ be an algebraically closed field with arbitrary characteristic. Then every unimodal isolated complete intersection singularity (ICIS) in $F[[x,y]]$ has the form in Table \ref{table-ord2}, \ref{table-ord2-char2}, \ref{table-ord3}, \ref{table-ord3-char2}, \ref{table-ord3-char3}. Besides Table \ref{table-ord2}, \ref{table-ord3-char2}, every form in the other tables is unimodal.  If additionally Conjecture \ref{conjecture} (resp. Conjecture \ref{conjecture2}) holds, then all the unimodal ICIS in Table \ref{table-ord2} (resp. Table \ref{table-ord3-char2}) are presented in Table \ref{table-ord2-mod1} (resp. Table \ref{table-ord3-char2-mod1}).
\end{theorem}

\newcommand{\etalchar}[1]{$^{#1}$}


\begin{thebibliography}{GPB{\etalchar{+}}08}

\bibitem[Arn76]{arnold-class-C}
V.I. Arnold.
\newblock Local normal forms of functions.
\newblock {\em Invent. Math.}, 35:87--109, 1976.


\bibitem[AVGZ12]{arnold}
V.I. Arnold, A.N. Varchenko, and S.M. Gusein-Zade.
\newblock  Singularities of differentiable maps. Volume 1. 
\newblock Reprint of the 1985 edition. {\em Modern Birkhäuser Classics.} Birkhäuser/Springer, New York, 2012. xii+382 pp. ISBN: 978-0-8176-8339-9. 

\bibitem[BGM12]{InvariantsOH}
Y.~Boubakri, G.M. Greuel, and T.~Markwig.
\newblock Invariants of hypersurface singularities in positive characteristic.
\newblock {\em Rev. Mat. Complut.}, 25 (2012), no. 1, 61–85. 

\bibitem[BKdP99]{complete-trans}
J.~Bruce, N.~Kirk, and A.~du~Plessis.
\newblock Complete transversals and the classification of singularities.
\newblock {\em Nonlinearity}, 10 (1997), no. 1,  253–275. 

\bibitem[DG83]{complete-unimodal-plane}
A.~Dimca and C.G. Gibson.
\newblock Contact unimodular germs from the plane to the plane.
\newblock {\em The Quarterly Journal of Mathematics}, 34(2)  (1983), 281--295.

\bibitem[Giu83]{germany-sing}
M.~Giusti.
\newblock Classification des singularit\'es isol\'ees simples d'intersections compl\`etes.
\newblock In {\em Singularities, {P}art 1 ({A}rcata, {C}alif., 1981)}, volume~40 of {\em Proc. Sympos. Pure Math.}, pages 457--494. Amer. Math. Soc., Providence, RI, 1983.


\bibitem[GN16]{right-simple}
G.M. Greuel and H.D. Nguyen.
\newblock Right simple singularities in positive characteristic.
\newblock {\em J. Reine Angew. Math.}, 712 (2016), 81–106.


\bibitem[GPB{\etalchar{+}}08]{singular-introduction}
G.M. Greuel, G.~Pfister.
\newblock  A singular introduction to commutative algebra.
\newblock Second, extended edition. {\em Springer, Berlin,} 2008. xx+689 pp. ISBN: 978-3-540-73541-0. 

\bibitem[Mat68]{mather-c-mapping}
J.N.~Mather.
\newblock Stability of $C^\infty $ mappings, {III.} {Finitely} determined map-germs.
\newblock {\em Publications Math\'ematiques de l'IH\'ES}, 35:127--156, 1968.

\bibitem[Ngu13]{phdclassification}
H.D. Nguyen.
\newblock  Classification of singularities in positive characteristic.
\newblock PhD thesis, 04 2013.


\bibitem[PG19]{finitedeter}
T.H. Pham and G.M. Greuel.
\newblock On finite determinacy for matrices of power series.
\newblock {\em Math. Z.} 290 (2018), no. 3-4, 759–774. 


\bibitem[PPG25]{icissimple}
T.H. Pham, G.~Pfister, and G.M. Greuel.
\newblock Classification of simple 0-dimensional isolated complete intersection singularities.
\newblock {\em Int. Math. Res. Not.}, 14 (2025), 1-22.

\bibitem[WR05]{actions-algebraic-group}
F. S.~Walter and A.~Rittatore.
\newblock  Actions and invariants of algebraic groups.
\newblock {\em Pure and Applied Mathematics (Boca Raton)},269. Chapman \& Hall/CRC, Boca Raton, FL, 2005. xvi+454 pp. ISBN: 978-0-8247-5896-7; 0-8247-5896-X.

\end{thebibliography}
\end{document}